\documentstyle[12pt]{amsart}
\newcommand{\nc}{\newcommand}

\nc{\C}{{\cal C}}
\renewcommand{\O}{{\cal O}}
\nc{\St}{\operatorname{St}^{\bullet}}
\nc{\B}{{\cal B}}
\nc{\A}{{\frak A}}
\nc{\D}{{\frak D}}

\nc{\N}{{\mathbf N}}
\nc{\CC}{{\mathbf C}}
\nc{\R}{{\mathbf R}}
\nc{\k}{{\mathbf  k}}
\nc{\Q}{{\mathbf Q}}
\nc{\U}{{\mathbf U}}
\nc{\Z}{{\mathbf Z}}
\nc{\Rhom}{\operatorname{RHom}\bul}
\nc{\Ad}{\operatorname{Ad}}
\nc{\Res}{\operatorname{Res}}
\nc{\gr}{\operatorname{gr}}
\nc{\tr}{\operatorname{tr}}
\nc{\End}{\operatorname{End}}
\nc{\r}{{\operatorname{re}}}
\renewcommand{\c}{{\operatorname{cont}}}
\renewcommand{\i}{{\operatorname{im}}}
\nc{\g}{{\frak g}}
\nc{\hatg}{\hat{\frak g}}

\nc{\gl}{{\frak g\frak l}}
\nc{\n}{{\frak n}}
\nc{\si}{{\frac\infty 2}}
\nc{\p}{{\frak p}}
\nc{\h}{{\frak h}}
\nc{\gothu}{{\frak u}}
\nc{\Ind}{\operatorname{Ind}}
\nc{\ch}{\operatorname{ch}}
\nc{\Coind}{\operatorname{Coind}}
\nc{\opp}{{\operatorname{opp}}}
\nc{\Ker}{\operatorname{Ker}}
\nc{\im}{\operatorname{Im}}
\nc{\Coker}{\operatorname{Coker}}
\nc{\dirlim}{\underset{\rightarrow}{\operatorname{lim}}}
\nc{\invlim}{\underset{\leftarrow}{\operatorname{lim}}}

\nc{\Ext}{\operatorname{Ext}^{\bullet}}
\nc{\ext}{\operatorname{Ext}}
\nc{\tilW}{\til{W}}

\nc{\lth}{\ell t}
\nc{\BB}{\cal{B}}
\nc{\Tor}{\operatorname{Tor}_{\bullet}}
\nc{\tor}{\operatorname{Tor}}
\nc{\Tors}{\operatorname{Tor}_{\frac \infty 2+\bullet}}
\nc{\Exts}{\operatorname{Ext}^{\frac \infty 2+\bullet}}
\nc{\Hom}{\operatorname{Hom}^{\bullet}}
\nc{\ad}{\operatorname{ad}}
\renewcommand{\hom}{\operatorname{Hom}}

\renewcommand{\mod}{\operatorname{-mod}}
\nc{\Mod}{\operatorname{Mod}}

\nc{\Barb}{\operatorname{Bar}^{\bullet}}

\nc{\upX}{X^{\uparrow}}
\nc{\upcD}{{\cal D}^{\uparrow}}
\nc{\upD}{D^{\uparrow}}
\nc{\dX}{X^{\downarrow}}
\nc{\dcD}{{\cal D}^{\downarrow}}
\nc{\dD}{D^{\downarrow}}
\nc{\upC}{{\cal C}^{\uparrow}}
\nc{\dC}{{\cal C}^{\downarrow}}
\nc{\M}{{\cal M}}
\nc{\underA}{\underline{A}}
\nc{\underB}{\underline{B}}
\nc{\underk}{\underline{\k}}
\nc{\Db}{D^{\bullet}}
\nc{\tenl}{\overset{\operatorname{L}}\ten}
\nc{\map}{\longrightarrow}

\nc{\bs}{\bigskip\\}
\nc{\ms}{\smallskip\\}

\nc{\tilBar}{\widetilde{\operatorname{Bar}}}
\nc{\tilBarb}{\widetilde{\operatorname{Bar}}^{\bullet}}
\nc{\linBar}{\overline{\operatorname{Bar}}}
\nc{\overr}{\overline{R}}
\nc{\overI}{\overline{I}}
\nc{\overX}{\overline{X}}
\nc{\overh}{\overline{h}}
\nc{\overY}{\overline{Y}}
\nc{\overW}{\overline{W}}
\nc{\linBarb}{\overline{\operatorname{Bar}}^{\bullet}}

\nc{\til}{\widetilde}
\nc{\oppA}{A^{\sharp}}

\nc{\Lemma}{{\bf Lemma:\ }}
\nc{\Theorem}{{\bf Theorem:}\ }
\nc{\Cor}{{\bf Corollary:}\ }
\nc{\Def}{{\bf Definition:}\ }
\nc{\Prop}{{\bf Proposition:}\ }
\nc{\Con}{{\bf Conjecture:}\ }
\nc{\Rem}{{\bf Remark:}\ }
\nc{\dok}{{\bf Proof.}\ }

\nc{\bul}{^{\bullet}}
\nc{\stand}{C^{\frac \infty 2+\bullet }}
\nc{\Cone}{\operatorname{Cone}}
\nc{\supp}{\operatorname{supp}}
\nc{\ten}{\otimes}

\nc{\ssn}{\subsection{}}
\nc{\sssn}{\subsubsection{}}

\nc{\hgt}{\operatorname{ht}}
\nc{\sqbinom}{\fracwithdelims[][0pt]}

\address{Independent University of Moscow, Pervomajskaya st. 16-18,
Moscow 105037, Russia}
\email{hippie@@ium.ips.ras.ru}

\thanks{Partially supported by Soros foundation.}

\author{S.~M.~Arkhipov}
\title{Semiinfinite cohomology of quantum groups II}
\begin{document}
\maketitle
\section{Introduction}
In [FF] B. Feigin and E. Frenkel introduced semiinfinite analogues of
the classical Bernstein-Gelfand-Gelfand (BGG) resolutions of
integrable simple modules over affine Lie algebras. These resolutions
are two sided complexes consisting of direct sums of so called
Wakimoto modules  suitable for computation of semi-infinite
homology of infinitely twisted Borel subalgebras in affine Lie
algebras (see e.~g.~[F] for the definition of the Lie algebra
semiinfinite homology). Feigin and Frenkel also suggested to consider
Wakimoto modules as direct limits of so called twisted Verma modules.
In [Ar3] it was shown that the semi-infinite BGG resolution itself can
be constructed as a direct limit of {\em twisted BGG resolutions} of
the integrable  simple module and that twisted BGG resolutions are
obtained from the classical BGG resolution with the help of
the {\em twisting functors}.

\ssn
The present paper is devoted to the analogue of this construction for
affine quantum groups.  Given a simply connected root datum
$(X,Y,\ldots)$ of the affine type $(I,\cdot)$ we consider an associative
algebra $\U$ over the field of rational functions $\Q(v)$ introduced
by Drinfeld and Jimbo and called the affine quantum group. This
algebra is a quantization of the universal enveloping algebra of the
affine Lie algebra $\hatg$ corresponding to $(I,\cdot)$. The algebra
$\U$ has  a natural {\em triangular decomposition} $\U=\U^-\ten
\U^0\ten\U^+$ as a vector space, where $\U^-$ (resp. $\U^0$, resp.
$\U^+$) denotes the quantum analogue of the universal enveloping
algebra for the standard negative nilpotent subalgebra (resp. for the
Cartan subalgebra, resp. for the standard positive nilpotent
subalgebra) in $\hatg$.  The infinitely twisted nilpotent  subalgebra
in $\hatg$ has a natural quantum analogue $\U^{\si-}$ in $\U$ as well.
Integrable simple modules $L(\lambda)$ also have nice quantizations
over $\U$. We preserve the notation $L(\lambda)$ for these
$\U$-modules.  The aim of this paper is to calculate semiinfinite
homology (introduced in the associative algebra case in [Ar1], [Ar2])
of $\U^{\si-}$ with coefficients in $L(\lambda)$.

\ssn
Let us describe briefly the structure of the paper. In the second
section we recall standard combinatorial definitions and constructions
concerning affine root systems and affine Weyl groups.  In particular
we define the length function on the affine Weyl group, the twisted
length function on $W$ with the twist $w\in W$, the semi-infinite
length function, the Bruhat order on $W$ and the semi-infinite Bruhat
order.  

In the third section we collect standard facts about affine quantum 
groups following mainly [L1] and [BeK]. We recall the construction of 
the automorphisms $T_w$, $w\in W$, of the affine quantum group $\U$ 
generating the action of the affine braid group $\B$ corresponding to 
$W$. We discuss $\theta_y$-twisted quantum nilpotent subalgebras  
$T_{\theta_y}(\U^-)$ in the affine quantum group $\U$ for an element 
of the transvection lattice $\theta_y\in T\subset W$ and their 
triangular decompositions following [Be], [BeK].

In the fourth and the fifth sections we collect nessesary definitions
and facts about  semiinfinite cohomology of a graded associative
algebra $A$ equipped with a triangular decomposition datum
(see~\ref{setup}) following [Ar1] and [Ar2]. In particular we
investigate the behaviour of semiinfinite homology and cohomology with
respect to the change  of the algebra $A$ (see
Propositions~\ref{negext},~\ref{neg},~\ref{posext},~\ref{pos}). 

In the sixth section we consider semiregular modules over the affine 
quantum group $\U$  with respect to various subalgebras in $\U$ and
endomorphism algebras of these modules. We prove that the endomorphism
algebra of the semiregular $\U$-module with respect to the subalgebra
$\U_{\theta_y}^-:=\U^-\cap T_{\theta_y}^{-1}(\U^+)$ contains $\U$ as a
subalgebra (see Theorem~\ref{bimodule}). Thus we obtain semiregular
$\U$-bimodules $S_{\U_{\theta_y}^-}^{\U}$ enumerated by elements of
the transvection lattice. Like in the affine Lie algebra case we
define twisting functors with the help of these bimodules. 

In the seventh section we recall quantum BGG resolutions 
$B\bul(\lambda)$ of the simple $\U$-modules $L(\lambda)$ following 
[M]. Like in [Ar3] for affine Lie algebras, we define quantum twisted 
BGG resoluitons as images of the complexes $B\bul(\lambda)$ under the 
twisting functors.  Using quantum twisted BGG resolutions we show that 
the semiinfinite homology spaces of the algebras $T_{\theta_y}(\U^-)$ 
with coefficients in $L(\lambda)$ are the same as in the affine Lie 
algebra case (see Lemma~\ref{answer}) and have a natural base 
enumerated by elements of $W$ graded by $\theta_y^{-1}$-twisted 
length. 

In the eighth section we show that the semiinfinite homology spaces of 
the algebras $T_{\theta_{mx}}(\U^-)$ ($x\in Y^{\prime\prime+}\subset 
T\subset W$ is a generic element, and  $m$ tends to infinity) with 
coefficients in $L(\lambda)$ form a projective system with the limit 
equal to the semiinfinite homology of the infinitely twisted quantum 
affine nilpotent algebra $\U^{\si-}$ with coefficients in 
$L(\lambda)$. In particular we prove that the latter semiinfinite 
homology spaces  are the same as in the affine Lie algebra case and 
have a natural base enumerated by elements of $W$ graded by the 
semi-infinite length function (see Theorem~\ref{stabil} and 
Corollary~\ref{final}).  Unfortunately we do not know the construction 
of the projective system on the level of quantum twisted BGG 
resolutions that would provide a quantum analogue of the semi-infinite 
BGG resolution.

\section{Notations and combinatorics of affine Weyl groups}
In this section we recall basic terminology and notations  concerning root systems,
affine Weyl groups following mainly [L1].

\subsection{Root data.} A {\em Cartan datum} is a pair $(I,\cdot)$ consisting
of a finite set $I$ and a symmetric bilinear form $\nu,\nu'\mapsto\nu\cdot\nu'$ on the
free abelian group $\Z[I]$ with values in $\Z$ such that $i\cdot i\in\{2,4,6,\ldots\}$
for any $i\in I$ and the number $2i\cdot j/i\cdot i\in\{0,-1,-2,\ldots\}$ for any
$i\ne j\in I$.

A Cartan datum $(I,\cdot)$ is said to be {\em of finite type} if the symmetric matrix
$(i\cdot j)$ indexed by $I\times I$ is positive definite. A Cartan datum $(I,\cdot)$
is said to be {\em of affine type} if it is irreducible and the symmetric
matrix $(i\cdot j)$ indexed by $I\times I$ is positive semi-definite but
not positive definite.

\sssn
A {\em root datum} of type $(I,\cdot)$ consists, by definition,
of two finitely generated abelian groups
$Y, X$ with a perfect bilinear pairing $\langle,\rangle:\ Y\times X\map\Z$ and
a pair of imbeddings $I\subset X\ (i\mapsto i')$ and $I\subset Y\
(i\mapsto i)$ such that $\langle i,j'\rangle =2i\cdot j/i\cdot i$ for all
$i,j\in I$.  We say that the root datum $(Y,X,\ldots)$ is {\em simply
connected} if $Y=\Z[I]$ with the obvious imbedding $I\map Y$;
$X=\hom(Y,\Z)$ with the obvious pairing $\langle,\rangle;\ Y\times X\map\Z$
and with the imbedding $I\map X$ defined by $\langle i,j'\rangle=2i\cdot j/i\cdot i$.

\sssn                           \label{fix}
We fix a simply connected root datum $(X,Y,\ldots)$ of the affine type $(I,\cdot)$.
We set $d_i=i\cdot i/2$. We suppose that the chosen affine Cartan datum is untwisted,
i.~e. there exists $i_0\in I$ such that $d_{i_0}=1$. Let $\overI:=I\setminus\{i_{0}\}$.
It is known that $(\overI,\cdot)$ is a Cartan datum of the finite type.
Moreover the chosen  root datum $(X,Y,\ldots)$ restricts to a root
datum $(\overX,\overY,\ldots)$ of the finite type $(\overI,\cdot)$.
Let $D=\max_id_i$, we have $D\in\{1,2,3\}$ and for each $i$,  $d_i$ is equal
either to $1$ or $D$. We define $\hat{d}_i$ by $d_i\hat{d}_i=D$ for  all
$i\in I$.
There are uniquely defined strictly positive integers $r_i,r'_i,\ i\in I$,
such that

(i) $\sum_ir_i\langle i,j'\rangle=0$ for all $j$ and $r_{i_0}=1$,\\
(ii) $\sum_jr'_j\langle i,j'\rangle=0$ for all $i$ and $r'_{i_0}=1$.

The dual Coxeter number corresponding to the Cartan datum $(I,\cdot)$ is defined as
follows: $h^{\vee}:=\underset{i\in I\setminus\{i_0\}}{\sum}r'_i$.

Let $\{\omega_i|i\in I\}$
be the basis of $X$ dual to $I\subset Y$. Consider an element
$c:=\sum_I r_ii\in Y$.  Then we have $\langle c,j'\rangle=0$ for all
$j\in I$ and $\sum_ir'_ii'=0$.

\subsection{Affine Weyl groups.}
For $i\in I$ we define reflections
$$
s_i:\ Y\map Y,\ s_i(y):=y-\langle y,i'\rangle i,\text{ and }
s_i:\ X\map X,\ s_i(x):=x-\langle i,x\rangle i'.
$$
The subgroup $W\subset \operatorname{Aut}Y$ generated by the reflections
$s_i,\ i\in I$, is called the {\em affine Weyl group}. We identify $W$
with the subgroup
in $\operatorname{Aut}(X)$ generated by $s_i,\ i\in I$. The subgroup
$\overline{W}\subset W$ generated by the reflections $s_i,\
i\in\overI$, is called the (finite) Weyl group corresponding to the
Cartan datum $(\overI,\cdot)$.

\sssn
Consider the set $R_\r$ (resp. $\overline{R}$) of elements of $Y$ of the form
$w(i)$ for some $i\in I$ and $w\in W$
(resp. of the form
$w(i)$ for some $i\in\overI$ and $w\in\overline{W}$).
It is called  the set of {\em real affine roots}
   (resp. the {\em finite root system}).
Let $R'$ (resp. $\overline{R}'$) be the set of vectors of $X$ of the form
$\omega(i')$ for some $i\in I$ (resp. for some $i\in\overI$ and
some $\omega\in \overW$). The assignment $i\mapsto i'$ extends
uniquely to a map
$$
\alpha\mapsto \alpha',\ R_\r\map R'\text{, such that  } \omega(y)'=\omega(y'),\
\omega\in W,\ y\in Y.
$$
The map restricts to bijection of $\overline{R}$ to  $\overline{R}'$.

There is a unique function $\alpha\mapsto d_{\alpha}$ on $R_\r$ such that it is
$W$-invariant and $d_i$ is defined in~\ref{fix}. We define $\hat{d}_{\alpha}$ by
$d_{\alpha}\hat{d}_{\alpha}=D$ for all $\alpha\in R_\r$.
Then we have

(i) $R_\r=\{\alpha+\hat{d}_{\alpha}mc|\alpha\in\overline{R},m\in\Z\}$,
 $R_\r+c=R_\r$;\\
(ii) $R'=\overline{R}'$,
$(\alpha+\hat{d}_{\alpha}mc)'=\alpha'$ for all $\alpha\in
\overline{R},\ m\in\Z$.

Consider also the set of {\em imaginary affine roots} $R_\i:=\{mc|m\in\Z\setminus\{0\}\}$.
The {\em affine root system} $R$ of the type $(I,\cdot)$ is defined as a union of
$R_\r$ and $R_\i$.  Note that $W$ acts trivially on $R_\i$.

\sssn
For $\alpha\in R_\r$ we denote by $s_\alpha$ the element of $W$  given by the reflection
in $Y$ (resp. in $X$)
$$
s_\alpha(y)=y-\langle y,\alpha'\rangle \alpha\text{ (resp. }
s_\alpha(x)=x-\langle \alpha,x\rangle \alpha').
$$
For any $\alpha\in\overline{R}$ and $m\in\Z$ we set $s_{\alpha,m}=s_h\in W$,
where $h=\alpha+\hat{d}_{\alpha}mc$.
Let
$Y'\subset X$ be a free abelian group
generated by the set
$\{i'|i\in \overI\}$ and let
$Y''\subset Y'$ be
a free abelian group
generated by the set
$\{\hat{d}_ii'|i\in\overI\}$. For $z\in Y''$ consider  a
transvection $\theta_z:\ X\map X$ given by $\theta_z(x)=x+\langle
c,x\rangle z$.

\sssn
\Lemma
For $\alpha\in\overline{R}_\r$ and $m\in\Z$ we have $s_{\alpha,0}\circ
s_{\alpha,m} =\theta_{\hat{d}_{\alpha}m\alpha'}$.  \qed

In particular $\theta_z\in W$ for any $z\in Y''$. Consider the map of the sets
$
\theta:\ Y''\map W,\
z\mapsto \theta_z,
$
and denote its image by $T\subset W$.
Then it is known that
the map $\theta$ is an injective homomorphism of groups,
$T$ is a normal subgroup in $W$ and
$W$ is a semidirect product of $T$ and $\overline{W}$. Note also that
$\theta_{\hat{d}_{\alpha}m\alpha'}$ acts on $Y$ by
$$
\theta_{\hat{d}_{\alpha}m\alpha'}(y)=y-\langle y,\alpha'\rangle\hat{d}_{\alpha}
mc.
$$
\sssn
As usual define the weight $\rho\in X$ by
$\rho(i)=1$ for all $i\in I$. Consider the {\em dot} action of
$W$ on $X$:
$$
w\cdot\lambda=w(\lambda+\rho)-\rho,\ w\in W,\ \lambda\in X.
$$
Then  the dot action of the Weyl group preserves the sets
$X_k:=\{\lambda\in X|\langle c,\lambda\rangle=k\}$.
From now on we denote the set of dominant weights at the level $k$
by $X^+_k$:
$$
X^+_k:=\{\lambda\in X_k|\langle i,\lambda\rangle\ge0,\ i\in I\}.
$$
\sssn
Recall that  the length of an element of the affine Weyl group $w$ is defined
as follows:
$$
\lth(w):=\sharp\{\alpha\in R^+|w^{-1}(\alpha)\in R^-\}.
$$
\Rem
The length of $w\in W$ is equal also to the minimal possible length of
expression of $w$ via the generators $s_i,\ i\in I$ (the length of a
reduced expression).

\sssn            \label{length}
For $w\in W$ consider a finite set
$R_w^-:=\{\alpha\in R^-|w(\alpha
)\in R^+\},\ R_w^+:=-R_w^-$. Let $w_1$ and $w_2$ be elements of the Weyl group
such that  $\lth(w_1)+\lth(w_2)=\lth(w_2w_1)$. Then $R_{w_2w_1}^{\pm}=
R_{w_2}^{\pm}\sqcup R_{w_1}^{\pm}$. Thus $R_{w_2}^{\pm}\cap w_1(R_{w_1}^{\mp})=\emptyset$, where the
set $-S\subset R$ consists of elements opposite to the ones of $S\subset R$.
Recall that we have $\underset{\alpha\in R_w^+}{\sum}\alpha=\rho-w^{-1}(\rho).$

\subsubsection{Extended affine Weyl group and affine braid group.}
We will need some extension of $W$.  Let $\Omega$ be the group of automorphisms of
$(W, I)$ whose restriction to $T$ is a conjugation by some element of $W$. Then
it is known that $\Omega$ is a finite group in correspondence with a certain group
of automorphisms of the Dynkin graph of $(I,\cdot)$ (see e.~g.~[BeK]). Then the {\em extended} affine
Weyl group $\tilW$ is defined as a semidirect product of $W$ and $\Omega$.
It is known that the decomposition of $W$ into the semidirect product of $T$ and $\overW$
can be extended to a decomposition of $\tilW$ into a semidirect product of
$\overX$ and $\overW$.

The length function on $W$ extends to the one on $\tilW$ by setting
$\lth(\tau w)=\lth(w)$ for $\tau\in \Omega$.

The {\em braid group} $\BB$ associated to $W$ is the group on generators
$T_w,\ w\in W$ with the relation $T_wT_{w'}=T_{ww'}$ if $\lth(ww')=
\lth(w)+\lth(w')$. The {\em extended braid group} $\til{\BB}$ is associated
to $\tilW$ in a similar way (see [Be]).

\subsection{Convex order on positive roots.}
Let $Y^+$ (resp. $Y^-$) be the set of all elements in $Y$  such that all their
coefficients with respect to the basis $I\subset Y$ are nonnegative
(resp. nonpositive).  Define the set of positive (resp. negative) roots
by $R^+=R\cap Y^+$, (resp. by $R^-=R\cap Y^-$).  We have
$R_\r^+=\{\alpha+\hat{d}_{\alpha}mc|\alpha\in\overline{R},m>0\}\sqcup
\overline{R}^+$ and $R_\i^+=\{mc|m\in\Z_{>0}\}$.

The height of a positive root $\alpha=\underset{i\in I}{\sum}b_ii$
is defined as follows: $\hgt\alpha=\underset{i\in I}{\sum}b_i$.
We extend the height function on the whole $Y$ by linearity.

\sssn
Let $w=s_{i_1}\ldots s_{i_k}$ be the reduced expression of an element of
the affine Weyl group $w\in W$. Then it is known that the set
$
\{i_k,s_{i_k}(i_{k-1}),\ldots,s_{i_k}s_{i_{k-1}}\ldots s_{i_2}(i_1)\}
$
coincides with the set $\{\alpha\in R^+|w(\alpha)\in R^-\}$.

\label{conv}
Recall the following construction of the set $R_\r^+$ (see [Be], [Pa]).
Fix an element $\theta_x\in T\subset W$ such that $\langle i,x\rangle>0$ for all
$i\in I$ and fix a reduced expression of $\theta_x$ via the generators of the affine Weyl group
$\theta_x=s_{j_1}\ldots s_{j_d}$.
Then in particular the opposite element of $\theta_x$ in $W$ has a reduced expression
$\theta_{-x}=s_{j_d}\ldots s_{j_1}$.
Let $(p_k)_{k\in\Z}$ be the sequence of integers
such that $p_k=j_{(k\operatorname{mod}(d))}$. Define the positive roots $\beta_k$ by
$$
\beta_k=s_{p_0}s_{p_{-1}}\ldots s_{p_{k+1}}(i_{p_k})\text{ for }k\le0\text{, and }
\beta_k=s_{p_1}s_{p_{2}}\ldots s_{p_{k-1}}(i_{p_k})\text{ for }k>0.
$$
\sssn
\Lemma

\qquad(i) The roots $\beta_k$ are distinct and the set
$\{\beta_k\}_{k\in\Z}$ coincides with $R_\r^+$;

\qquad(ii) each subsection $s_{p_k}\ldots s_{p_l}$ for $k<l$ is reduced.\qed

A total order on the set of positive roots is called {\em convex} if for  any
$\alpha\in R^+$ and $\beta\in R^+$ we have
$
\alpha<\alpha+\beta<\beta\text{ in the order if }\alpha+\beta\in R^+.
$

\sssn
\Lemma
The order on $R^+$ defined by
$$
\beta_0<\beta_{-1}<\beta_{-2}<\ldots<rc<\ldots<sc<\ldots<\beta_2<\beta_1
$$
is convex. Here $R_\r^+$ is identified with the set $\{\beta_k\}_{k\in\Z}$
using the previous Lemma, and $r<s$.\qed

\subsection{Semi-infinite Bruhat order on the Weyl group.}
We say that $w'$ follows $w$ in the Weyl group if there exist a reduced
expression of $w'$ and  $p\in\{1,\ldots,\lth w'\}$ such that
$$
w'=s_{i_1}\ldots s_{i_{\lth w'}},\ w=s_{i_1}\ldots s_{i_{p-1}}s_{i_{p+1}}\ldots
s_{i_{\lth w'}}
$$
and $\lth w=\lth w'-1$.

Recall that the usual Bruhat order on the Weyl group is the partial order on
$W$ generated by the relation ``$w'$ follows $w$ in $W$''. It is denoted by
$\ge$.

\sssn
\Lemma
The relation
$
\{$there exists $\lambda_0\in Y^{\prime\prime+},\langle i,\lambda_0\rangle>0$
 for every
$i\in\overline{I}$, such that for any
$\lambda\in Y^{\prime\prime+}$ we have
$\theta_{\lambda}\theta_{\lambda_0}w'
\ge
\theta_{\lambda}\theta_{\lambda_0}w\}
$
is a partial order on $W$.
\qed

\sssn
\Def
We call this partial order the semi-infinite Bruhat order on $W$ and denote it
by $\ge^{\si}$.

\sssn
\Rem
The semi-infinite Bruhat order defined above in fact coincides with the
partial order on the affine Weyl group defined in [L2], Section 3, in terms
of combinatorics of {\em alcoves} in $Y'\ten\R$. However we will not need
the comparison statement. Further details on the partial order can be
found in [L2].

\sssn  \label{tlength}
Denote the set $\{\alpha\in R|\alpha=\beta+\hat{d}_{\beta}mc,\
\beta\in\overr^+,m\in\Z\}$ (resp.  the set $\{\alpha\in
R|\alpha=\beta+\hat{d}_{\beta}mc,\ \beta\in\overr^-,m\in\Z\}$) by $R^{\frac
\infty 2+}_\r$ (resp.  by $R^{\frac \infty 2-}_\r$).
By definition we set $R^{\si+}:=R^{\si+}_\r\sqcup R^+_\i$ and
$R^{\si-}:=R^{\si-}_\r\sqcup R^-_\i$.
Following [FF] we introduce
the semi-infinite length function on the affine Weyl group as follows:
$$
\lth^{\si}(w):=
\sharp\{\alpha\in R^{\si+}\cap R^+|w(\alpha)\in R^-\}
-
\sharp\{\alpha\in R^{\si-}\cap R^+|w(\alpha)\in R^-\}.
$$
Next we consider  the twisted length function on the affine Weyl group.
For $u,w\in W$ set
$$
\lth^w(u):=\lth(w^{-1}u)-\lth(w^{-1}).
$$
In particular for $\mu\in -Y^{\prime\prime+}$ and $\nu\in Y^{\prime\prime+}$
we have $\lth^{\theta_{\mu}}(\theta_{\nu})=\lth(\theta_{\nu})$.

\sssn
\Lemma
\label{maincomb} For every $w_1,w_2\in W$ there exists $\mu_0\in
-Y^{\prime\prime+}$, $\langle i,\mu_0\rangle<0$ for every $i\in\overline{I}$,
such that for every $\mu\in -Y^{\prime\prime+}$ we have
$$
\lth(\theta_{-\mu}\theta_{-\mu_0}w_1)-
\lth(\theta_{-\mu}\theta_{-\mu_0}w_2)=
\lth^{\si}(w_1)-
\lth^{\si}(w_2).
$$
In particular for every $\mu\in -Y^{\prime\prime+}$ we have
$\lth^{\theta_{\mu}\theta_{\mu_0}}(w)=\lth^{\si}(w)$.
\qed

\sssn
\Cor \label{qq}
For every $x\in Y^{\prime\prime+},\ \langle i,x\rangle>0$ for every $i\in \overline{I}$,
for every $w\in W$ there exists $m_0\in \N$ such that for every $m>m_0$ we have
$\lth^{\theta_{mx}^{-1}}(w)=\lth^{\si}(w)$.

\dok
Fix $\mu_0=\mu_0(w)$ from the previous Lemma. There exists $m_0\in\N$ such that for every $m>m_0$
we have $\mu:=-mx-\mu_0\in-Y^{\prime\prime+}$. Thus by the previous Lemma we obtain
$$
\lth^{\theta^{-1}_{mx}}(w)=\lth^{\theta_{-mx}}(w)=\lth^{\theta_{\mu}\theta_{\mu_0}}(w)
=\lth^{\si}(w).\qed
$$
By \ref{tlength} and the previous Corollary for $\mu\in Y^{\prime\prime+}$ we have
$\lth^{\si}(\theta_{\mu})=\lth(\theta_{\mu})$.

\section{Affine quantum groups.}
In this section we present some facts about triangular decompositions of
various subalgebras in affine quantum groups following mainly [L1] and [BeK].

\ssn
Assume that a root datum $(Y,X,\ldots)$ of the affine type $(I,\cdot)$ is given.
We consider the associative $\Q(v)$ algebra $'\til{\U}$ with $1$ defined by generators
$$
E_i\quad(i\in I),\quad F_i\quad(i\in I),\quad K_\mu\quad(\mu\in Y)
$$
and the following relations:
\begin{gather*}
K_0=1,\ K_\mu K_{\mu'}=K_{\mu+\mu'}\text{ for all }\mu,\mu'\in Y;\\
K_\mu E_i=v^{\langle\mu,i'\rangle}E_iK_\mu\text{ for all }i\in I,\ \mu\in
Y;\ K_\mu F_i=v^{-\langle\mu,i'\rangle}F_iK_\mu\text{ for all }i\in I,\
\mu\in Y;\\
E_iF_j-F_jE_i=\delta_{ij}\frac{\til{K}_i-\til{K}_{-i}}{v_i-v_i^{-1}}.
\end{gather*}
Here for $\nu=\sum_i\nu_ii\in Y$ we set $\til{K}_{\nu}:=\prod_iK_{(i\cdot i/2)\nu_i}$
and $v_i:=v^{i\cdot i/2}$.

Thus $'\til{\U}$ is $Y$-graded. We denote the $Y$-grading of a homogenious element $u\in{}'\til{\U}$
by $\deg^Y u\in Y$. We define the quantum  adjoint  action
of $E_i$ and $F_i$ on a homogenious element $u\in{}'\til{\U}$ as follows:
$$
\Ad_{E_i}^q(u):=E_iu-v^{\langle\deg^Y u,i'\rangle}uE_i\text{ and }
\Ad_{F_i}^q(u):=F_iu-v^{-\langle\deg^Y u,i'\rangle}uF_i.
$$
\sssn
\Def
The associative algebra obtained from $'\til{\U}$ by taking quotient by the ideal generated
by {\em quantum Serre relations}
$$
\Ad_{E_i}^{q\ -\langle i,j'\rangle+1}(E_j)\text{ and }
\Ad_{F_i}^{q\ -\langle i,j'\rangle+1}(F_j)\text{ for all } i,j\in I
$$
is called the {\em affine quantum group} $\til{\U}$.

We denote  by $\til{\U}^-$ (resp. by $\til{\U}^+$, resp. by $\til{\U}^0$) the subalgebra in
$\til{\U}$ generated by $F_i$ for all $i\in I$ (resp. by $E_i$ for all $i\in I$,
resp.  by all $K_\mu$, $\mu\in Y$).  It is known that the multiplication
in $\til{\U}$ provides the isomorphisms of vector spaces
$\til{\U}^-\ten\til{\U}^0\ten\til{\U}^+\til{\map}\til{\U}$ and $\til{\U}^+\ten\til{\U}^0\ten\til{\U}^-\til{\map}\til{\U}$
(see~e.~g.~[L1]). By definition set $\til{\U}^{\ge0}:=\til{\U}^0\ten\til{\U}^+$ and
$\til{\U}^{\le0}:=\til{\U}^0\ten\til{\U}^-$. Evidently $\til{\U}^{\ge0}$ and $\til{\U}^{\le0}$ are
subalgebras in $\til{\U}$.

\sssn \label{inv}
We introduce the $\Q$-algebra antiautomorphism $\kappa$ of $\til{\U}$ defined by:
$$
\kappa(E_i)=F_i,\ \kappa(F_i)=E_i,\ \kappa(K_\nu)=K_{-\nu},\ \kappa(v)=v^{-1},
\ i\in I,\ \nu\in Y.
$$
\subsection{Convex PBW bases in $\til{\U}$.}
Now we recall the definition of the analogues of Lusztig generators of $\til{\U}^+$
introduced in the affine Cartan datum case by Beck (see [Be]).
First we construct the real root generators.

\sssn
\Lemma (see e.~g. [BeK], 1.5)
The formulas
$$
T_iE_i=-F_iK_i,\ T_i(E_j)=\frac{(-1)^{\langle i,j'\rangle}}{[\langle i,j'\rangle]_{d_i}!}
(\Ad_{E_i}^q)^{-\langle i,j'\rangle}(E_j)\text{ if }i\ne j,\ T_iK_\nu=K_{s_i\nu},\
i,j\in I,\nu\in Y,
$$
define an action of the affine braid group $\BB$ on $\til{\U}$ by automorphisms that
can be extended to the action of $\til{\BB}$.\qed

Fix a total convex order on $R^+$ constructed starting from a positive
transvection $\theta_x\in T$ (see~\ref{conv}). For each
$\beta_k\in R_\r^+$ define the root vector $E_{\beta_k}\in \til{\U}_{\beta_k}$
as follows:
$$
E_{\beta_k}=T_{i_0}^{-1}\ldots T_{i_{k+1}}^{-1}(E_{i_k})\text{ if }k\le0,\
E_{\beta_k}=T_{i_1}\ldots T_{i_k}(F_{i_k})\text{ if }k>0.
$$
We will need also negative real root generators. For $\beta_k\in R_\r^+$ define the root vector
$F_{\beta_k}\in\til{\U}_{\beta_k}$ by
$$
F_{\beta_k}=T_{i_0}^{-1}\ldots T_{i_{k+1}}^{-1}(F_{i_k})\text{ if }k\le0,\
F_{\beta_k}=T_{i_1}\ldots T_{i_k}(E_{i_k})\text{ if }k>0.
$$
\sssn
Next we construct the imaginary root generators.
For $i\in \overline{I}$ and $m>0$ set
$$
\psi_m^{(i)}:=\til{K}_i^{-1}[E_i,E_{mc-i}],\
\psi_{-m}^{(i)}:=\kappa(\psi_m^{(i)}),\ \psi_0^{(i)}=\frac{\til{K}_i-\til{K}_i^{-1}}
{v_i-v_i^{-1}}.
$$
\sssn
For $k>0$ and $i\in \overline{I}$ define imaginary root vectors $E_{kc}^{(i)}$ by the following
functional equation:
$$
\exp\left((v_i-v_i^{-1})\sum_{k=1}^{\infty}E_{kc}^{(i)}t^k\right)=
1+(v_i-v_i^{-1})\sum_{k=1}^{\infty}\psi_k^{(i)}t^k.
$$
\sssn
\Rem
As defined $E_{\beta_k}\in\til{\U}^{\ge0}$  and $F_{\beta_k}\in \til{\U}^{\le0}$
for $\beta_k\in R_\r^+$. Set
$\dot{E}_{\beta}:=K_{-\alpha}E_{\beta}$ if $\beta=-\alpha+k\hat{d}_{\alpha}c$,
and $\alpha\in\overr^-$, otherwize set $\dot{E}_{\beta}:=E_\beta$ (still $\beta\in R^+_\r$).
The negative real root generators $\dot{F}_{\beta}$ are obtained from $F_{\beta}$ in
a similar way.

Consider a set $\til{R}^+:=\left(R^{\si-}_\r\cap R^+\right)\sqcup\underbrace
{R^+_\i\sqcup\ldots\sqcup R^+_\i}_{\sharp\overline{I}}\sqcup \left(R^{\si+}_\r\cap R^+\right)$.
Denote the obvious projection $\til{R}^+\map R^+$ by $\pi$. Fix an
arbitrary total order on $\til{R}^+$ such that $\pi$ would become a
map of the totally ordered sets. Note that $\til{R}^+$ enumerates the
positive root vectors in $\til{\U}^+$ constructed above. The sets
$\til{R}^-$, $\til{R}^{\si+}$ and $\til{R}^{\si-}$ are defined in a similar
way.

Abusing notation we write $\langle\alpha,\beta'\rangle$ for $\langle\pi(\alpha),
\pi(\beta)'\rangle$, where $\alpha,\beta\in \til{R}^+$.

\sssn
\Lemma (see [Be])        \label{basis}
The elements $\dot{E}_{\beta}$, $\beta\in \til{R}^+$, belong to $\til{\U}^+$ and generate it.
Let $\Z_{\ge0}^{\til{R}^+}$ be the set of all $\Z_{\ge0}$-valued functions on $\til{R}^+$ with finite
support. For $(a_\beta)\in\Z_{\ge0}^{\til{R}^+}$ consider a monomial $M_{(a_\beta)}\in\til{\U}^+$
as folows: $M_{(a_\beta)}:=\underset{\beta\in \til{R}^+}
{\prod}\dot{E}_{\beta}^{a_{\beta}}$ (taken in the chosen convex order
on the positive roots). Then such monomials form a base of $\til{\U}^+$ over
$\Q(v)$.\qed

By analogy with the Lie algebra case we call the base $\{M_{(a_\beta)}\}$
the convex PBW base of $\til{\U}^+$.  The main purpouse of introducing convex
PBW bases is the following statement.

\sssn
\Theorem  (see [Be])                  \label{multiplication}
Let $\alpha,\beta\in \til{R}^+$ be such that $\beta>\alpha$ in the total convex order
on $\til{R}^+$. Then we have
$$
\dot{E}_{\beta}\dot{E}_{\alpha}-v^{\langle\alpha,\beta'\rangle}
\dot{E}_{\alpha}\dot{E}_{\beta}=\underset{\alpha<\gamma_1<\ldots
<\gamma_n<\beta}{\sum}a_{\gamma}\dot{E}_{\gamma_1}^{b_1}\ldots
\dot{E}_{\gamma_n}^{b_n},
$$
where $\gamma$ denotes the vector $(\gamma_1,\ldots,\gamma_n)$, $\gamma_s\in \til{R}^+,$
and the coefficients $a_{\gamma}\in \Q[v,v^{-1}]$.\qed

A similar statement holds for $\til{\U}^-$.

\subsubsection{Filtration on $\til{\U}^+$.} Using the PBW type base we define a filtration on
$\til{\U}^+$ like in [DCK]. The notion of a $S$-filtration for a totally ordered set $S$ was introduced there.
We view $\Z_{\ge0}^{\til{R}^+}$ as a totally ordered semigroup with the usual lexicographical order.
Introduce a $\Z_{\ge0}^{\til{R}^+}$-filtration on $\til{\U}^+$ by letting $F^s\til{\U}^+$
be the span of the monomials
$$
M_{(a_{\beta})}\in\til{\U}^+,\ (a_{\beta})=(\ldots,0,\ldots,0,a_{\beta_m},\ldots,a_{\beta_n},0\ldots,0,\ldots)\in \Z_{\ge0}^{\til{R}^+}
$$
such that $(a_{\beta})<s$ in the total order.

\sssn
\Prop
(see [BeK]) \label{graded}
The associated $\Z_{\ge0}^{\til{R}^+}$-graded algebra $\gr\til{\U}^+$ of the
$\Z_{\ge0}^{\til{R}^+}$-filtered algebra $\til{\U}^+$ is an algebra over $\Q(v)$ on generators
$E_\alpha$, $\alpha\in \til{R}^+$, and relations $E_\alpha E_\beta =v^{\langle\alpha,\beta'\rangle}E_\beta E_\alpha$
for $\beta<\alpha$ in the total convex order.\qed

\subsection{Subalgebras in $\til{\U}^+$.}
Fix a positive root $\beta_k\in R_\r^+$.
We denote the set
$\{\alpha|\alpha\in \til{R}^+,\ \alpha<\beta_k\}$ (resp.
$\{\alpha|\alpha\in \til{R}^+,\ \alpha\ge\beta_k\}$) by $\til{R}^+_{<\beta_k}$
(resp. by $\til{R}^+_{\ge\beta_k}$).
Let $\til{\U}^+_{<\beta_k}$ (resp.
$\til{\U}^+_{\ge\beta_k}$) be the subalgebra in $\til{\U}^+$ generated by the root elements
$\{\dot{E}_\alpha|\alpha\in \til{R}^+_{<\beta_k}\}$
(resp. by the root elements
$\{\dot{E}_\alpha|\alpha\in \til{R}^+_{\ge\beta_k}\}$).
The corresponding subalgebras $\til{\U}^-_{<-\beta_k}$ and $\til{\U}^-_{\ge-\beta_k}$
in $\til{\U}^-$ are defined in  a similar way.
Then we have the following statement.

\sssn
\Lemma            \label{ddd}
The monomials
$M_{(a_\beta)}:=\underset{\beta\in \til{R}^+_{<\beta_k}}
{\prod}\dot{E}_{\beta}^{a_{\beta}}$
(resp. the monomials
$M_{(a_\beta)}:=\underset{\beta\in \til{R}^+_{\ge\beta_k}}
{\prod}\dot{E}_{\beta}^{a_{\beta}}$)
taken in the chosen convex order
on the positive roots form a base of $\til{\U}^+_{<\beta_k}$
(resp. of $\til{\U}^+_{\ge\beta_k}$) over
$\Q(v)$. The algebra $\til{\U}^+$ admits a decomposition $\til{\U}^+\til{\map}\til{\U}^+_{<\beta_k}
\ten\til{\U}^+_{\ge\beta_k}$ as a $\Q(v)$-vector space.

\dok
The Lemma follows immediately from Lemma~\ref{basis} and
Theorem~\ref{multiplication}.\qed

In particular consider the positive roots $\beta_{md}=
s_{p_1}\ldots s_{p_{md-1}}(i_{p_{md}})$ and
$\beta_{-md}=s_{p_0}s_{p{-1}}\ldots s_{p_{-md+1}}(i_{p_{md}})$ (see~\ref{conv}).
Then we have
$$
\til{R}^+_{\ge\beta_{md}}=\{\alpha\in \til{R}^+|\theta_{-mx}(\pi(\alpha))\in R^-\}
\text{ and }
\til{R}^+_{<\beta_{-md-1}}=\{\alpha\in \til{R}^+|\theta_{mx}(\pi(\alpha))\in R^-\}.
$$
We denote the subalgerba $\til{\U}^+_{<\beta_{-md-1}}$ (resp. the subalgebra
$\til{\U}^-_{<-\beta_{-md-1}}$) by $\til{\U}^+_{\theta_{mx}}$ (resp. by $\til{\U}^-_{\theta_{mx}}$.
Consider the antiautomorphism
$
T_{\theta_{mx}}^{-1}=(T_{p_d}^{-1}\ldots T_{p_1}^{-1})^m
$
of the algebra $\til{\U}$.

\sssn                                  \label{im}
\Lemma (see [Be], 3.10, Proposition 2)
We have $T_{\theta_x}(\dot{E}_{mc}^{(i)})=\dot{E}_{mc}^{(i)}$ for all $i\in \overline{I}$
and $m\in\Z\setminus \{0\}$.\qed

\sssn
\Lemma
We have  $T_{\theta_{mx}}^{-1}(\til{\U}^+_{\ge\beta_{md}})=\til{\U}^-_{\theta_{mx}}$ and
$T_{\theta_{mx}}^{-1}(\til{\U}^-_{\ge-\beta_{md}})=\til{\U}^+_{\theta_{mx}}$. In particular
$\til{\U}^-_{\theta_{mx}}=\til{\U}^-\cap T^{-1}_{\theta_{mx}}(\til{\U}^+)$
and
$\til{\U}^+_{\theta_{mx}}=\til{\U}^+\cap T^{-1}_{\theta_{mx}}(\til{\U}^-)$.

\dok
The first two statements are proved by the following calculation.
Let $\dot{E}_{\beta_l}$ (resp. $\dot{F}_{\beta_l}$), $0<l\le md$, be
one of the root generators of $\til{\U}^+_{\ge\beta_{md}}$ (resp. of $U^-_{\ge-\beta_{md}}$).
Then, forgetting aboul the factor from $\til{\U}^0$, we have
\begin{gather*}
T_{\theta_{mx}}^{-1}(E_{\beta_l})=
T_{p_{md}}^{-1}\ldots T_{p_1}^{-1}(T_{p_1}\ldots T_{p_l}(F_{p_l}))=
T_{p_{md}}^{-1}\ldots T_{p_{l+1}}^{-1}(F_{p_l})\\ =
T_{p_{0}}^{-1}\ldots T_{p_{l-md+1}}^{-1}(F_{p_{l-md}})=F_{-\beta_{l-md}}.
\end{gather*}
A similar calculation holds for $F_{\beta_l}$. Thus we see that the subalgebras
$T_{\theta_{mx}}^{-1}(\til{\U}^+_{\ge\beta_{md}})$ and $\U^-_{\theta_{mx}}$
(resp.
$T_{\theta_{mx}}^{-1}(\til{\U}^+_{\ge\beta_{md}})$ and $\U^-_{\theta_{mx}}$)
have the same generators.
The last two statements follow from the previous considerations. \qed

\sssn
\Lemma  \label{dec}
We have
$\til{\U}^+_{\ge\beta_{-md-1}}=\til{\U}^+\cap T^{-1}_{\theta_{mx}}(\til{\U}^+)$.

\dok
The fact that the real root generators of the two subalgebras coincide is proved
the same way as in the previous Lemma. The imaginary root generators are invariant under the automorphism
$T_{\theta_{mx}}^{-1}$ by Lemma~\ref{im}. Now use Lemma~\ref{basis}. \qed

In particular we obtain a triangular decomposition of the algebra
$T_{\theta_{mx}}^{-1}(\til{\U}^+)$
(resp.
$T_{\theta_{mx}}^{-1}(\til{\U}^-)$):
\begin{gather*}
T_{\theta_{mx}}^{-1}(\til{\U}^+)
\til{\map}\til{\U}^+_{\ge\beta_{-md-1}}\ten
\til{\U}^-_{\theta_{mx}}=
T_{\theta_{mx}}^{-1}(\til{\U}^+)\cap\U^+\ten
\til{\U}^-_{\theta_{mx}}\\
\text{(resp. }
T_{\theta_{mx}}^{-1}(\til{\U}^-)
\til{\map}\til{\U}^-_{<\beta_{md}}\ten
\til{\U}^+_{\theta_{mx}}=
T_{\theta_{mx}}^{-1}(\til{\U}^-)\cap\U^-\ten
\til{\U}^+_{\theta_{mx}}).
\end{gather*}
as a $\Q(v)$-vector space.

Consider the $\Z_{\ge0}^{\til{R}^+}$-filtration on the algebra $\til{\U}^+_{\theta_{mx}}$
obtained by restriction from the one on $\til{\U}^+$. Note that by definition
nontrivial graded quotient spaces are enumerated in fact by the totally
lexicographically ordered {\em finitely generated} semigroup $\Z_{\ge0}^{R_{\theta_{mx}}^+}\subset
\Z_{\ge0}^{\til{R}^+}$.

\sssn
\Lemma \label{assgr}
The associated $\Z_{\ge0}^{R_{\theta_{mx}}^+}$-graded algebra $\gr\til{\U}^+_{\theta_{mx}}$ of the
$\Z_{\ge0}^{R_{\theta_{mx}}^+}$-filtered algebra $\til{\U}^+_{\theta_{mx}}$ is an algebra over $\Q(v)$ on generators
$E_\alpha$, $\alpha\in R_{\theta_{mx}}$ and relations $E_\alpha E_\beta
=v^{\langle\alpha,\beta'\rangle}E_\beta E_\alpha$ for $\beta<\alpha$ in the
total convex order.\qed

\ssn \label{level}
It is easily checked that the element $\til{Z}:=\underset{i\in I}{\prod}\til{K}_i^{r_i}$
is central in $\til{\U}$. To simplify the exposition we add a certain root
of the element $\til{Z}$ to the algebra $\til{\U}$. Namely we set
$\U:=\til{\U}[Z]$ where $Z$ is defined as  a central element such that $Z^{Dh^{\vee}}=\til{Z}$.
Again we have a triangular decomposition of the algebra $\U$: $\U=\U^-\ten\U^0\ten\U^+$
as a $\Q(v)$-vector space where as before $\U^+$ (resp. $\U^-$, resp.
$\U^0$) denotes the subalgebra in $\U$ generated by $E_i,\ i\in I$ (resp.
by $F_i,\ i\in I$, resp. by $K_\mu,\ \mu\in Y$, and $Z$).

Fix the level $k\in\Z$. We define the algebra $\U_k$
as the  quotient algebra of $\U$ by the relation $Z=v^k$.

Note that $\U^{\pm}=\til{\U}^{\pm}$.
All the statements about convex bases in $\til{\U}$ from the previous subsection
hold for $\U$.

\section{Semiinfinite cohomology of associative algebras.}
In this section we present some definitions and statements concerning
semiinfinite cohomology of associative algebras following mainly [Ar2],
sections 2 and 3.

\ssn
\label{setup}%
Suppose we have a graded associative algebra
$A=\underset{n\in \Z}{\bigoplus}A_n$ over a field $\k$.
Let $B$ and $N$ be graded
subalgebras in $A$ satisfying the following conditions:

 \qquad(i) $N$ is positively graded;

 \qquad(ii) $N_0=\k;$

 \qquad(iii) $\dim N_n<\infty$ for any $n\in\N;$

 \qquad(iv) $B$ is negatively graded;

 \qquad(v) the multiplication in $A$ defines the isomorphisms of graded
 vector spaces
$$
 B\otimes N\map A\text{ and }N\otimes B\map A.
$$
In particular $N$ is naturally augmented.
We denote the augmentation ideal
$\underset{n>0}{\bigoplus}N_n$ by $\overline{N}$.

\sssn
The category of left graded $A$-modules with morphisms that
preserve gradings is denoted by $A\mod$. We define the
functor of the grading shift
$$
 A\mod\map A\mod:\ M\mapsto
 M\langle i\rangle,\ M\langle i\rangle_m:=M_{i+m}, i\in\Z.
$$
The space $\underset
{i\in\Z}{\bigoplus}\hom_{A\mod}(M_1,M_2\langle i\rangle)$ is
denoted by $\hom_A(M_1,M_2)$.

\sssn
We fix a left $B$-augmentation on $A$ provided by the
isomorphism of graded left $B$-modules $B\cong
A\ten_N\underk$ where $\underk:=N/\overline{N}$ is the
trivial $N$-module. The $A$-module $A\ten_N\underk$ is denoted by
$\underB$.

\sssn
We introduce certain subcategories in the category of
complexes ${\cal K}om(A\mod)$.  For
$M\bul\in{\cal K}om(A\mod)$ the support of $M$ is defined as
follows:
$$
 \supp M\bul:=\{(p,q)\in\Z^{\oplus 2}|M_p^q\ne
 0\}.
$$
For $s_1,s_2,t_1,t_2\in \Z,\ s_1,s_2>0$, the set $
\{(p,q)\in \Z^{\oplus 2}|s_1q+p\ge t_1,\ s_2q-p\ge t_2\} $
(resp. the set $ \{(p,q)\in \Z^{\oplus 2}|s_1q+p\le t_1,\
s_2q-p\le t_2\}) $ is denoted by $\upX(s_1,s_2,t_1,t_2)$
(resp. by $\dX(s_1,s_2,t_1,t_2)$).

Let $\upC(A)$ (resp. $\dC(A)$) be the full subcategory in
${\cal K}om(A\mod)$ consisting of complexes $M\bul$ that
satisfy the following condition:

\quad(U) there exist $s_1,s_2,t_1,t_2\in \Z,\ s_1,s_2>0$,
such that $\supp M\bul\subset\upX(s_1,s_2,t_1,t_2)$ (resp.

\quad(D) there exist $s_1,s_2,t_1,t_2\in \Z,\ s_1,s_2>0$,
such that $\supp M\bul\subset\dX(s_1,s_2,t_1,t_2)$).

\subsection{Relative bar resolutions.}
The standard bar resolution $\tilBarb (A,B,M)\in {\cal
K}om(A\mod)$ of an $A$-module $M$
with respect to the subalgebra $B$ is by definition the following complex
of $A$-modules
$$
  \tilBar^{-n}(A,B,M):=A\otimes_B\ldots \otimes_BA\otimes_BM \
  (n+1 \text{ times}),
$$
with the standard bar differential.
Consider the subspace $\linBarb (A,B,M)$ in $\tilBarb(A,B,M)$ as follows:
$$
  (\linBar)^{-n}(A,B,M):=\{ a_0\otimes\ldots\otimes a_n\otimes
  v\in \tilBar^{-n}(A,B,M)|\ \exists \ s\in\{1,\ldots,n\}:\
  a_s\in B\}.
$$
Evidently it is a submodule in $\tilBarb (A,B,M)$ preserved by the
differential.

The quotient complex $\Barb (A,B,M):=\tilBarb (A,B,M)/\linBarb
(A,B,M)$ is called the restricted bar resolution of the
$A$-module $M$ with respect to the subalgebra $B$.
For a complex of $A$-modules $M\bul\in\upC(A)$ its relative
restricted bar resolution is defined as a total complex of
the bicomplex $\Barb (A,B,M\bul)$.

Let $M\in \upC(A)$. Then  it is known that the restricted relative bar
resolution of $M$ also belongs to $\upC(A)$ and is quasiisomorphic to $M$
(see [Ar2], Lemma 2.3.4).

\subsection{Relative cobar DG-algebra.}
Here we present a construction of a canonical DG-algebra representing
$\Rhom_A(\underB,\underB)$.
Consider the relative restricted bar resolution $\Barb
(A,B,\underB)$ of the $A$-module $\underB$ and the complex
of graded vector spaces
$$
 \Db(A,B):=\Hom_A(\Barb(A,B,\underB),\underB).
$$
Clearly $\Db(A,B)\cong
\Hom_\k(\underset{n\ge0}{\bigoplus}\overline{N}^{\ten
n},\underB)$ as a vector space, and $D(A,B)$ belongs to
$\dC({\cal V}ect)$.

\sssn \label{mult}%
We introduce a structure of a DG-algebra on $\Db(A,B)$.
First we define a DG-algebra structure on
$\Hom_A(\tilBarb(A,B,\underB),\underB)$. Note that by
Schapiro lemma
$$
  \Hom_A(\tilBar^{-m}(A,B,\underB),\underB)=
  \Hom_B(A\ten_B\ldots\ten_BA\ten_B\underB,\underB)\ (m\text{ times}).
$$
Let $f\in \hom^m_A(\tilBarb(A,B,\underB),\underB),\ g\in
\hom^n_A(\tilBarb(A,B,\underB),\underB)$, i.~e.
$$
  f:\
  \underbrace{A\ten_B\ldots\ten_BA}_m\ten_B\underB\map
  \underB,\ g:\
  \underbrace{A\ten_B\ldots\ten_BA}_n\ten_B\underB\map
  \underB,
$$
both $f$ and $g$ are $B$-linear. By definition
set $f\cdot g:\
\underbrace{A\ten_B\ldots\ten_BA}_{m+n}\ten_B\underB\map
  \underB:$
$$
  (f\cdot g)(a_1\ten\ldots\ten a_{m+n}\ten b):=
  f(a_1\ten\ldots\ten a_n\ten g(a_{n+1}\ten\ldots\ten
  a_{m+n}\ten b)).
$$
Then the multiplication
equips $\Hom_A(\tilBarb(A,B,\underB),\underB)$ with a
DG-algebra structure and the subcomplex
$$
   \Db(A,B):=\{f\in\Hom_A(\tilBarb(A,B,\underB),\underB)\ |
  f\equiv 0 \text{ on }
  \linBarb(A,B,\underB)\subset\tilBarb(A,B,\underB)\}
$$
is a  DG-subalgebra in
$\Hom_A(\tilBarb(A,B,\underB),\underB)$ (see [Ar2], Lemma 2.4.2).

\sssn
The algebra $\Db(A,B)$ has a kind of a trianuglar decomposition. First
its zero grading component
$D^0(A,B)=\hom_A(A\ten_B\underB,\underB)=
\hom_B(\underB,\underB)\cong B^{\opp}$ as an algebra (yet
the inclusion $B^{\opp}\hookrightarrow \Db(A,B)$ is not a
morphism of DG-algebras --- the differential in $\Db(A,B)$
does not preserve $B^{\opp}$).

Next consider the induction functor $\Ind_N^A:\ N\mod\map
A\mod$. Then the canonical map
\begin{multline*}
  \Hom_N(\Barb(N,\k,\underk),\underk)\map\Hom_A(\Ind_N^A(\Barb(N,\k,\underk)),
  \Ind_N^A(\underk))\\
  \cong\Hom_A(\Barb(A,B,\underB),\underB)=\Db(A,B)
\end{multline*}
is an inclusion of DG-algebras.
Denote its image by $\Db(N,\k)$. Then we have
$
\Db(A,B)=\Db(N,\k)\ten B^{\opp}
$
as a graded vector space.

\sssn \label{setup2}
We will need another restriction on the algebra $A$ as follows.

(vi) Both maps $\Db(N,\k)\ten B^{\opp}\map\Db(A,B)$
and $B^{\opp}\ten\Db(N,\k)\map\Db(A,B)$ provided by the multiplication
in $\Db(A,B)$ are isomorphisms of vector spaces.

\subsubsection{Free algebras over $B$.}
\label{free}
Fix a bimodule $M$ over the algebra $B$. Denote by $T_B(M)$ the {\em free}
algebra over $B$ generated by $M$. By definition $T_B(M)$ is the
algebra on generators space equal to $M$ and with relations as
follows:
\begin{gather*}
\left(m_1\ten\ldots\ten
m_k\right)\cdot\left(m_{k+1}\ten\ldots\ten m_{k+l}\right)= m_1\ten\ldots\ten
m_{k+l},\\
b\cdot m=bm,\ m\cdot b=mb;\ m,m_1,\ldots,m_{k+l},bm,mb\in
M, b\in B.
\end{gather*}
Thus we have
$
T_B(M)=\underset{n\ge0}{\bigoplus}\underbrace{M\ten_B\ldots\ten_B
M}_{n}
$
as a vector space.
Consider the $B$-bimodules $A/B$ and
$M:=\hom_B(A/B,B)$. The latter one is provided by two $B$-module
structures on the space of homomorphisms between the  {\em left}
$B$-modules:  the left $B$-module structure is given by the {\em
right} $B$-multiplication in $A/B$ and the right $B$-module structure
is given by the right $B$-multiplication in $B$. Then evidently we
have an isomorphism of algebras $\Db(A,B)\cong T_B(M)$.

\subsection{Bar duality  functors.}
\label{functor}%
Next we construct canonical $\Db(A,B)$-DG-modules representing
$\Rhom_A(\underB,*)$ and $*\tenl_A\underB$.
Let  $M\bul\in\dC(A),\ M^{\prime\bullet}\in\upC(A^\opp)$.
By definition set
$$
  \dD(A,B,M\bul):=\Hom_A(\Barb(A,B,\underB),M\bul),\
  \upD(A,B,M^{\prime\bullet}):=M^{\prime\bullet}\ten_A\Barb(A,B,\underB).
$$
Evidently the vector space $\dD(A,B,M\bul)$
 (resp. $\upD(A,B,M^{\prime\bullet})$) belongs to
$\dC({\cal V}ect)$ (resp. to $\upC({\cal V}ect)$).  Similarly
to~\ref{mult} we define the right action of $\Db(A,B)$ on
 $\dD(A,B,M\bul)$ (resp the left action of $\Db(A,B)$ on
$\upD(A,B,M^{\prime\bullet})$) (see [Ar2], 2.4.4).

Then the multiplication equips
$\Hom_A(\tilBarb(A,B,\underB),M\bul)$ (resp.
$M^{\prime\bullet}\ten_A\tilBarb(A,B,\underB)$)
with a structure of the right DG-module over $\Db(A,B)$
(resp. with a structure of the left DG-module over $\Db(A,B)$.
Note also that
$\dD(A,B,M\bul)
\subset\Hom_A(\tilBarb(A,B,\underB),M\bul)$
is a DG-submodule and
$\upD(A,B,M^{\prime\bullet})$ is a quotient
DG-module of $M^{\prime\bullet}\ten_A\tilBarb(A,B,\underB)$
(see [Ar2], Lemma 2.4.5).

\sssn
Denote the category of right (resp. left) DG-modules
$$
  M\bul=\underset{p,q\in\Z}{\bigoplus}M_p^q,\
  \operatorname{d}:\ M_p^q\map M_p^{q+1},
$$
over $\Db(A,B)$, with morphisms being
morphisms of DG-modules that preserve also the second
grading, by $\Db(A,B)\mod$ (resp. by $\Db(A,B)^\opp\mod$).
The subcategory in $\Db(A,B)\mod$ (resp. in
$\Db(A,B)^\opp\mod$) that consists of DG-modules satisfying
the condition (D) (resp. (U)) is denoted by $\dC(\Db(A,B))$
(resp.  by $\upC(\Db(A,B)^\opp)$.

The localizations of $\dC(A)$, $\upC(A^\opp)$,
$\dC(\Db(A,B))$, $\dC(\Db(A,B))$
and $\upC(\Db(A,B)^\opp)$ by the class of quasiisimorphisms
are denoted by $\dcD(A)$, $\upcD(A^\opp)$, $\dcD(\Db(A,B))$,
and $\upcD(\Db(A,B)^\opp)$ respectively.

By~\ref{functor} $\dD$ and $\upD$ define the  functors
$$
 \dD_A:\ \dC(A)\map \dC(\Db(A,B))\text{ and }\upD_A:\
 \upC(A^\opp)\map \upC(\Db(A,B)^\opp).
$$
\sssn
\Theorem  (see [Ar2], Theorem 2.4.7)
\label{equiv}

\qquad(i) The functor $\dD_A$ is well defined as a functor
from $\dcD(A)$ to $\dcD(\Db(A,B))$;

\qquad(ii) $\dD_A:\ \dcD(A)\map \dcD(\Db(A,B))$ is an
equivalence of triangulated categories;

\qquad(iii) the functor $\upD_A$ is well defined as a
functor from $\upcD(A^\opp)$ to $\upcD(\Db(A,B)^\opp)$;

\qquad(iv) $\upD_A:\ \upcD(A^\opp)\map \upcD(\Db(A,B)^\opp)$
is an equivalence of triangulated categories. \qed

\subsection{Construction of the algebra $\protect\oppA$.}
Now we are going to introduce an algebra $\oppA$ such that
$\oppA=B^{\opp}\ten N ^{\opp}$ as a vector space and
$\Db(A,B)^{\opp}\cong \Db(\oppA,B^{\opp})$.

Recall that the DG-algebra $\Db(A,B)$ is a tensor product of
the DG-subalgebra $\Db(N,\k)$ and the subalgebra (not a
DG-subalgebra) $B^{\opp}$. The isomorphism of vector spaces is provided
by the condition (vi) from \ref{setup2}.
Recall also that the algebra
$\Db(N,\k)$ is isomorphic to the tensor algebra over the
graded vector space $\overline{N}^*$ and the differential in
it is generated by the map
$\overline{N}^*\map\overline{N}^*\ten\overline{N}^*$ dual to
the multiplication map in $N$ and extended to the whole
algebra $T(\overline{N}^*)$ by the Leibnitz rule.

To define a DG-algebra $C^{\vee}$ such that $C^{\vee}$ is a
tensor product of its DG-subalgebra $T(\overline{N}^*)$ with
a differential given by the map dual to the multiplication
in $N$ and a (not DG-) subalgebra $B^{\opp}$ one has to
specify the following additional data:
\begin{itemize}
  \item a linear map $B\ten
  \overline{N}^*\map\overline{N}^*\ten B$ generating the
  multiplication in $C^{\vee}$; \item a linear map $B\map
  \overline{N}^*\ten B$ providing a component of the
  differential in $C^{\vee}$,
\end{itemize}
satisfying certain
constraints (that provide the associativity constraint and
the Leibnitz rule in the DG-algebra $C^{\vee}$).

On the other hand to define an algebra $C$ such that $C$ is
a tensor product of two its subalgebras $B$ and $N$ as a
vector space and the conditions~\ref{setup} are satisfied
one has to specify the following data (additional to the
algebra structures on $B$ and $N$):
\begin{itemize}
  \item a linear map $\overline{N}\ten B\map B\oplus
  B\ten\overline{N}$ providing the multiplication in $C$
\end{itemize}
satisfying certain constraints (that provide
the associativity constraint in the algebra $C$).

\sssn
\Prop
The construction of the dual algebra $C\mapsto
C^{\vee}:=\Db(C,B)$ provides a one to one correspondence
between the two described types of data, i.~e. for every
DG-algebra $C^{\vee}$ such that
$C^{\vee}=T(\overline{N}^*)\ten B^{\opp}$ as a vector space
there exists an algebra $C=B\ten N$ as a vector space such
that the DG-algebras $C^{\vee}$ and $\Db(C,B)$ are
isomorphic.
\qed

\sssn
\Lemma
The DG-algebras $\Db(N,\k)^\opp$ and $\Db(N^\opp,\k)$ are
isomorphic to each other.  \qed

Thus we have a triangular decomposition  for the DG-algebra
$\Db(A,B)^{\opp}$ as follows:
$$
\Db(A,B)^{\opp}=B\ten \Db(N^{\opp},\k)=\Db(N^{\opp},\k)\ten B.
$$
\sssn
\Cor
There exists an associative algebra $\oppA$ such that
$\oppA$ contains two subalgebras $B^{\opp}$ and $N^{\opp}$
satisfying  the conditions~\ref{setup} for $\oppA$, $B^\opp$
and $N^\opp$ such that the DG-algebra $\Db(A,B)^{\opp}\cong
\Db(\oppA,B^{\opp})$.  \qed

\sssn
\Cor
The functor $D_{\oppA}^\downarrow$ provides an equivalence
of triangulated categories \\
$\dcD(\oppA)\til{\map}\dcD(\Db(A,B)^\opp)$.  \qed

\subsection{Definitions of semiinfinite cohomology of associative
algebras.}
Consider a left graded $N$-module $N^*:=\underset{n\ge
0}{\bigoplus} \hom_\k(N_n,\k)$.  The action of $N$ on the
space is defined as follows.
$$
 f:\ N\map \k,\ n\in N,\text{
 then } (n\cdot f)(n'):=f(n'n).
$$
Consider also the left
$A$-module $S_N^A:=\Ind_N^A(N^*)=A\ten_NN^*$.  Evidently
$S_N^A\cong B\ten N^*$ as a $B$-module and by~\ref{setup}
$S_N^A\cong N^*\ten B$ as a $N$-module.

There is another left $A$-module with these two properties.
$S_N^{\prime A}:=\hom_B(A,B)$ with the left action of $A$ defined as
follows.
$$
 f:\ A\map B,\ a\in A,\text{ then }(a\cdot
 f)(a'):=f(a'a).
$$
\sssn
\Lemma
The $A$-modules $S_N^A$ and $S_N^{\prime A}$ are isomorphic.

\dok
Fix the decomposition $A=B\ten N$ provided by the
multiplication in $A$.  For $f\in N^*$ define $f_2:\ A\map
B$, $f_2(b\ten n):=f(n)b$.  Define the pairing
$$
 S_N^A\times A\map B,\ f\ten a\times a'\mapsto f_2(a'a).
$$
One checks directly the correctness of the definition.
Now the condition (vi) from \ref{setup2} provides that the defined map is
an isomorphism of vector spaces.
\qed

Thus $S_N^A\cong \Coind_B^AB$. The functors of induction from
$N$ to $A$  and of coinduction from $B$ to $A$ provide the
natural inclusions of algebras
$$
 N^{\opp}\hookrightarrow
 \hom_A(S_N^A,S_N^A)\text{ and }B^\opp\hookrightarrow
 \hom_A(S_N^A,S_N^A).
$$
Note that $\hom_A(S_N^A,S_N^A)=\hom_{\k}(N^*,B)$ as a vector
space. Below we construct an inclusion of associative algebras
$\oppA\subset\End_A(S_N^A)$ such that the subalgebras $B^{\opp}$
and $N^{\opp}$ in $\End_A(S_N^A)$ are exactly the ones provided by
the triangular decomposition of $\oppA$.

\sssn
Consider the $\Db(A,B)$-DG-module $\dD_A(S_N^A)$. By Theorem~\ref{equiv}
we have an isomorphism of associative algebras
$$
\End_A(S_N^A)\cong\ext_{\Db(A,B)}^0(\dD_A(S_N^A),\dD_A(S_N^A)).
$$
Here $\ext$ functors are taken in the derived category of right
$\Db(A,B)$-DG-modules, that is, in the category of {\em left}
$\Db(A,B)^{\opp}$-DG-modules. Consider also the left
$\Db(A,B)^{\opp}$-DG-module $\upD_{\oppA}(\underline{A}^{\sharp})$,
where $\underline{A}^{\sharp}$ denotes
the regular $\oppA$-module.

\sssn
\Lemma
$\upD_{\oppA}(\underline{A}^{\sharp})\cong\dD_A(S_N^A)$ in the derived category of right
$\Db(A,B)$-DG-modules.\qed

\sssn
\Prop  (see [Ar2], Lemma 3.4.5)
The functor $\upD_{\oppA}$ provides an {\em inclusion} of associative
algebras
$$
\oppA\map\ext_{\Db(A,B)^\opp}^0(\upD_{\oppA}(\underline{A}^{\sharp}),
\upD_{\oppA}(\underline{A}^{\sharp}))=\End_A(S_N^A).\qed
$$
Thus $S_N^A$ becomes a $A-\oppA$-bimodule.

\subsubsection{Continious $\hom$ description of $\oppA$.}
\label{contin}%
We intorduce  a topology
on $A$ (resp. on a graded $A$-module $M$) defined by the filtration $F^mA:=\underset{n<m}{\bigoplus}A_n$
(resp. by the filtration $F^mM:=\underset{n<m}{\bigoplus}M_n$). In particular the multiplication
in $A$ and the action of $A$ on $M$ are given by continious maps.
Denote the space of continious linear maps between two graded $A$-modules $M$ and $M'$
equipped with this topology by $\hom^\c(M,M')$. Thus we have
$$
\hom^\c(M,M')=\underset{n<0}{\bigoplus}\hom(M_n,M')\oplus\underset{n\ge0}{\prod}
\hom(M_n,M').
$$
For left graded $A$-modules $M$ and $M'$ consider the space of continious morphisms
$$
\hom^\c_A(M,M'):=\{f\in\hom^\c(M,M')|f(am)=af(m)\text{ for }a\in A,m\in M\}.
$$
In particular we have
\begin{gather*}
\hom^\c_A(S_N^A,S_N^A)=\hom^\c_B(S_N^A,B)=\hom^\c(N^*,B)\\=\underset{n\le0}{\bigoplus}
\hom((N^*)_n,B)=\underset{n\ge0}{\bigoplus}((N_n)^*,B)=N\ten B.
\end{gather*}
It is easily checked that the images of the inclusions $B^\opp\subset\End_A(S_N^A)$
and $N^\opp\subset\End_A(S_N^A)$ belong to the space of continious endomorphisms.
Thus we obtain the following statement.

\sssn
\Lemma
$\oppA=\hom^\c_A(S_N^A,S_N^A)$.\qed

Now we give a definition of associative algebra
semiinfinite cohomology amd compare it with the one
presented in [Ar1].

\sssn
\Def
Let $M\bul\in\upC(A^\opp)$,
$M^{\sharp\bullet}\in\dC(\oppA)$. Then  set
$$
  \Exts_A(M^{\sharp\bullet},M\bul):=
  \Ext_{\Db(A,B)^\opp}(\dD_{\oppA}(M^{\sharp\bullet}),\upD_A(M\bul)).
$$
Note that by definition the semiinfinite $\ext$ functor
maps complexes exact by either of the variables to zero.
We will need also the semiinfinite $\tor$ functor.

\sssn
\Def
Let $M\bul\in\dC(A)$,
$M^{\sharp\bullet}\in\dC(\oppA)$. Then  set
$$
  \Tors^A(M^{\sharp\bullet},M\bul):=
  \Tor^{\Db(A,B)}(\dD_{\oppA}(M^{\sharp\bullet}),\dD_A(M\bul)).
$$
Note that the definition can be rewritten as follows:
$$
  \Tors^A(M^{\sharp\bullet},M\bul):=
  \left(\Exts_A(M^{\sharp\bullet},M^{\bullet*})\right)^*.
$$
Here $M^{\bullet*}$ denotes the right $A$-module dual to $M\bul$.

The definition of semiinfinite cohomology in [Ar1] used the following
standard complex:
$$
  \stand(M^{\sharp\bullet},M\bul):=
  \Hom_{\oppA}\left(\Barb(\oppA,N^\opp,M^{\sharp\bullet}),
  \Barb(A^\opp,B^\opp,M\bul)\ten_AS_N^A\right).
$$
\sssn
\Theorem (see [Ar1], Theorem 3.6.2)
Let $M\bul\in\upC(A^\opp)$,
$M^{\sharp\bullet}\in\dC(\oppA)$. Then
$$
\Exts_A(M^{\sharp\bullet},M\bul)=H\bul(\stand(M^{\sharp\bullet},M\bul)).\qed
$$
We present also a similar statement for the semiinfinite $\tor$ functor.
Let $M\bul\in\dC(A)$, $M^{\sharp\bullet}\in\dC(\oppA)$.
We define the standard complex for the computation of the semiinfinite $\tor$
functor by
$$
C_{\si+\bullet}(M^{\sharp\bullet},M\bul):=
\left(C^{\si-\bullet}(M^{\sharp\bullet},M^{\bullet*})\right)^*.
$$
\sssn
\Cor
Let $M\bul\in\dC(A)$, $M^{\sharp\bullet}\in\dC(\oppA)$. Then
$$
\Tors^A(M^{\sharp\bullet},M\bul)=H^{-\bullet}(C_{\si+\bullet}(M^{\sharp\bullet},M\bul)).\qed
$$
\subsection{Choice of resolutions.}
Recall that an object $M\in A\mod$ is called {\em injective} (resp. {\em projective})
{\em relative to the
subalgebra} $N$ if for every complex of $A$-modules $X\bul$ such that
$X\bul$ is homotopic to zero as a complex of $N$-modules we have
$H\bul(\Hom_A(X\bul,M))=0$
(resp.
$H\bul(\Hom_A(M,X\bul))=0$).
An object $M\in A\mod$ is called {\em semijective}
if it is both $N$-projective and $A$-injective relative to $N$. Here we
recall that one can use semijective resolutions to calculate semiinfinite
cohomology.

\sssn
\Lemma (see [Ar1], Lemma 3.4.1) \label{extres}
The following facts hold
for $L\bul\in \dC (\oppA ),$ $M\bul\in
\upC (A^{\opp})$:

\qquad(i) if $M\bul$ is $N$-projective, then  we have
$$
\Exts_A(L\bul ,M\bul )=H\bul (\Hom_{\oppA }(\Barb(\oppA ,N,L\bul ),
M\bul\otimes_AS_N^A));
$$
\qquad(ii) if $L\bul$ is $B$-projective, then        we have
$$
\Exts_A(L\bul ,M\bul )=H\bul (\Hom_{\oppA }(L\bul ,
\Barb (A,B,M\bul )
\otimes_A
S_N^A));
$$
\qquad(iii) if $M\bul$ is semijective     then           we have
$$
\Exts_A(L\bul ,M\bul )=
H\bul (\Hom_{\oppA } (L\bul ,
M\bul\otimes_AS_N^A)).\qed
$$
We present also the analogue of Lemma~\ref{extres}(iii) for semiinfinite $\tor$
functors to be used later.
We call $M\bul\in\dC(A)$ {\em co-semijective} if it
is both $N$-injective and $A$-projective relative to $N$.

\sssn
\Cor    \label{torres}
Suppose $M\bul\in\dC(A)$ is co-semijective and
$L\bul\in\dC(\oppA)$.
Then we have
$$
\Tors^A(L\bul,M\bul)=H\bul(\Hom_A(S_N^A,M\bul)\ten_{\oppA}L\bul).\qed
$$
\sssn
\Prop (see [Ar1], Theorem 3.5)\label{derived}
Semiinfinite $\ext$ functor is well defined on the
corresponding derived categories:
$$
\Exts_A:\ {\cal
D}^{\downarrow}(\oppA )\times {\cal D}^{\uparrow}(A^{\opp})\map {\cal D}
({\cal V}ect).  \qed
$$
A similar fact holds evidently for semiinfinite $\tor$ functors.

\section{Functorial properties of semiinfinite cohomology}
Now we describe the behaviour of
semiinfinite cohomology with respect to a change of rings in spirit of
[Ar3], 6.2.

\subsection{Negative extension case.}
Suppose we have an inclusion of algebras $A\subset A'$ such that both $A$
and $A'$ have triangular decompositions
$A=B\ten N$ and $A'=B'\ten N$ as vector spaces satisfying the conditions
(i)-(v) from~\ref{setup} and (vi) from~\ref{setup2}. Suppose also that
$B$ is a subalgebra in $B'$ (and note that positive parts of the
triangular decompositions coincide). Denote the modules over the
corresponding algebras provided by the $B$-augmentation (resp. by the $B'$-augmentation)
by
$\underline{B}_A$ and $\underline{B}'_{A'}$ respectively. We have a natural
functor of induction $\Ind_A^{A'}:\ A\mod\map A'\mod$  mapping
$\underline{B}_{A}$ to $\underline{B}'_{A'}$.
Consider corresponding morphism of functors
$$
\Ind_A^{A'}:\ \Rhom_{A}(\underline{B}_{A},\underline{B}_{A})\map
\Rhom_{A'}(\underB'_{A'},\underB'_{A'}).
$$
It is represented by the morphism of DG-algebras $\Db(A,B)\map\Db(A',B')$
constucted as follows.

\sssn
\Lemma       \label{morph'}
We have $\Barb(A',B',\underB'_{A'})\cong\Ind_A^{A'}\Barb(A,B,\underB_A)$.
\qed

\sssn
\Cor
The natural map
\begin{multline*}
\Db(A,B,\underB)=\Hom_A(\Barb(A,B,\underB_A),\underB_A)\map \\
\Hom_{A'}(\Barb(A',B',\underB'_{A'}),\underB'_{A'})=\Db(A',B,\underB_{A'})
\end{multline*}
provided by the functor $\Ind_A^{A'}$ is a morphism of DG-algebras.\qed

\sssn \label{mor1'}
Suppose we have a right $A'$-module $M$.
By the previous Lemma we have an isomorphism of $\Db(A,B)$-DG-modules
\begin{gather*}
\upD_{A'}(M)=M\ten_{A'}\Barb(A',B',\underB'_{A'})\cong\\
M\ten_{A'}\Ind_A^{A'}\Barb(A,B,\underB_A)\cong
M\ten_A\Barb(A,B,\underB_A)=\upD_A(M).
\end{gather*}
\sssn                 \label{m1}
Next we construct the morphism of the $^\sharp$-algebras.
Consider  again the functor of induction
$$
\Ind_A^{A'}:\ A\mod\map A'\mod,\ M\mapsto A'\ten_AM.
$$
Then we have
$
S_N^{A'}=\Ind_N^{A'}N^*=\Ind_{A}^{A'}\Ind_N^AN^*=\Ind_A^{A'}S_N^A.
$
Thus we obtain an inclusion of algebras
$\End_A(S_N^A)\hookrightarrow\End_{A'}(S_N^{A'}).$

\sssn
\Lemma \label{canis1}
Let $X\bul\in\dC(B)$. Then we  have  acanonical isomorphism of $A'$-modules
$$
\Coind_{B'}^{A'}\circ\Ind_B^{B'}X\bul\cong\Ind_A^{A'}\circ\Coind_B^AX\bul.
$$
\dok
By~\ref{setup2}(vi) we have
$
\Coind_B^AX\bul\til{\map}\left(\Coind_B^A\underB_A\right)\ten_BX\bul\til{\map}S_N^A\ten_BX\bul.
$
Evidently, $\Ind_A^{A'}S_N^A=S_N^{A'}$. Thus we obtain
$
\Ind_A^{A'}\circ\Coind_B^AX\bul\til{\map}S_N^{A'}\ten_BX\bul.
$
On the other hand we have
$$
\Coind_{B'}^{A'}\circ\Ind_B^{B'}X\bul\til{\map}
\Coind_{B'}^{A'}\underB'_{A'}\ten_{B'}\Ind_B^{B'}X\bul=S_N^{A'}\ten_BX\bul.
$$
Taking the composition of the constructed isomorphisms we obtain the required one.
\qed

\sssn
\Lemma
The  map constructed in~\ref{m1} provides an inclusion of algebras
$\oppA\hookrightarrow A^{\prime\sharp}$ well defined with respect to the
triangular decompositions.

\dok
The inclusion $N^{\opp}\subset\End_A(S_N^A)$ (resp.
$N^{\opp}\subset\End_{A'}(S_N^{A'})$) is provided by the realisations of the
semiregular modules as coinduced ones from the coregular $N$-module. Thus
the constructed map $\End_A(S_N^A)\map\End_{A'}(S_N^{A'})$ identifies the
images of the algebra $N^{\opp}$.

Next by Lemma~\ref{canis1}  we have the  isomorphism of functors $\dC(B)\map
\dC(A')$:
$$
\Coind_{B'}^{A'}\circ\Ind_B^{B'}\cong\Ind_A^{A'}\circ\Coind_B^A.
$$
Thus the images of two inclusions
$$
B^{\opp}\hookrightarrow\End_A(S_N^A)\hookrightarrow\End_{A'}(S_N^{A'})
\text{ and }B^{\opp}\hookrightarrow B^{\prime\opp}
\hookrightarrow\End_{A'}(S_N^{A'})
$$
coincide. Now recall that as a vector space the subalgebra
$\oppA\subset\End_A(S_N^A)$ is described as $N^\opp\ten B^{\opp}$.\qed

Note also that the inclusions $A\subset A'$ and
$\oppA\subset A^{\prime\sharp}$ provide the same morphism of DG-algebras
$\Db(A,B)=\Db(A^{\sharp},B)^{\opp}\map\Db(A^{\prime\sharp},B')^{\opp}=\Db(A',B')$.

\sssn  \label{mor2'}
Suppose we have an $A^{\prime\sharp}$-module $M^{\sharp}$.
By Lemma~\ref{morph'} we have an isomorphism of $\Db(A,B)$-DG-modules
\begin{gather*}
\dD_{A^{\prime\sharp}}(M^{\sharp})=\Hom_{A^{\prime\sharp}}(\Barb(A^{\prime\sharp},B',
\underB'_{A^{\prime\sharp}}),M^{\sharp})\cong
\Hom_{A^{\prime\sharp}}(\Ind_{\oppA}^{A^{\prime\sharp}}\Barb(\oppA,B,
\underB_{\oppA}),M^{\sharp})\\ \cong
\Hom_{\oppA}(\Barb(\oppA,B,\underB_{\oppA}),M^{\sharp})=\dD_{\oppA}(M^{\sharp}).
\end{gather*}
\sssn
\Prop              \label{negext}
There exists a natural (functorial) morphism
$$
\Exts_{A'}(M^{\sharp},M)\map
\Exts_{A}(M^{\sharp},M).
$$
\dok
For any two $\Db(A',B')$-DG-modules $X\bul$ and $Y\bul$
the morphism of the DG-algebras $\Db(A,B)\map\Db(A',B')$ obtained from the morphism of functors
$\Rhom_{A}(*,*)\map\Rhom_{A'}(\Ind_A^{A'}(*),\Ind_A^{A'}(*))$ provides a natural
restriction map
$$
\Rhom_{\Db(A',B')}(X\bul,Y\bul)\map
\Rhom_{\Db(A,B)}(X\bul,Y\bul).
$$
In particular we have a natural morphism
$$
\Rhom_{\Db(A',B')}(\dD_{A^{\prime\sharp}}(M^{\sharp}),\upD_{A'}(M))\map
\Rhom_{\Db(A,B)}(\dD_{A^{\prime\sharp}}(M^{\sharp}),\upD_{A'}(M)).
$$
Now the required morphism is obtained by the composition of the constructed one with the isomorphisms
described in~\ref{mor1'} and~\ref{mor2'}.\qed

Similarly we obtain the following statement.
\sssn
\Prop       \label{neg}
Suppose we have an inclusion of algebras $A\subset A'$ such that both algebras have triangular decompositions
$A=B\ten N$ and $A'=B'\ten N$ as vector spaces and the inclusion preserves the triangular decompositions, i~e.
$B\subset B'$. Let $M\in\dC(A')$ and
$M^{\sharp}\in\dC(A^{\prime\sharp})$.  Then there exists a natural
(functorial) morphism
$$
\Tors^{A}(M^{\sharp},M)\map
\Tors^{A'}(M^{\sharp},M).\qed
$$
Note that by the definition of the semiinfinite $\tor$ functor we have
$\Tors^A(M^{\sharp},S_N^A)=M^{\sharp}$ for every $\oppA$-module $M^{\sharp}$.

\sssn
\Lemma       \label{canmap2}
let $M^{\sharp}\in\dC(A^{\prime\sharp})$. Then
the canonical map
$$
\Tors^{A}(M^{\sharp},S_{N}^{A'})\map\Tors^{A'}(M^{\sharp},S_{N}^{A'})
$$
coincides with the obvious one
$\Ind_{B^{\opp}}^{B^{\prime\opp}}\Res_{B^{\opp}}^{A^{\prime\sharp}}M^{\sharp}\map M^{\sharp}$.
\qed

\subsection{Positive extension case.}
Suppose we have an inclusion of algebras $A\subset A'$ such that both $A$
and $A'$ have triangular decompositions
$A=B\ten N$ and $A'=B\ten N'$ as vector spaces satisfying the conditions
(i)-(v) from~\ref{setup} and (vi) from~\ref{setup2}. Suppose also that
$N$ is a subalgebra in $N'$ (and note that negative parts of the
triangular decompositions coincide). Denote the modules over the
corresponding algebras provided by the $B$-augmentations by
$\underline{B}_A$ and $\underline{B}_{A'}$ respectively. We have a natural
functor of restriction $\Res_A^{A'}:\ A'\mod\map A\mod$  mapping
$\underline{B}_{A'}$ to $\underline{B}_A$.
Consider corresponding morphism of functors
$$
\Res_A^{A'}:\ \Rhom_{A'}(\underline{B}_{A'},\underline{B}_{A'})\map
\Rhom_A(\underline{B}_A,\underline{B}_A).
$$
It is represented by the morphism of DG-algebras $\Db(A',B)\map\Db(A,B)$
constucted as follows.

Consider the morphism of bar resolutions
$\Barb(A,B,\underB_A)\map\Barb(A',B,\underB_{A'})$ provided by the
inclusion of algebras.

\sssn
\Lemma       \label{morph}
The  composition of the described map with the restriction morphism
provides a morphism of DG-lagebras
\begin{gather*}
\Db(A',B)=\Hom_{A'}(\Barb(A',B,\underB_{A'}),\underB_{A'})\map
\Hom_A(\Barb(A',B,\underB_{A'}),\underB_A)\map\\
\Hom_A(\Barb(A,B,\underB_A),\underB_A) =\Db(A,B).\qed
\end{gather*}
Suppose we have a right $A'$-module $M$.
Consider the  morphism of functors
$$
\Res_A^{A'}:\ M\tenl_A\underline{B}_{A}\map
M\tenl_{A'}\underline{B}_{A'}.
$$
It is represented by a morphism of DG-modules over  $\Db(A',B)$
constucted similarly to the one in Lemma~\ref{morph}:
$\upD_A(M)\map\upD_{A'}(M)$. \label{mor1}

\sssn  \label{m2}
Consider also a the functor of coinduction
$$
\Coind_A^{A'}:\ A\mod\map A'\mod,\ M\mapsto\hom_A(A',M).
$$
Then we have
$
S_{N'}^{A'}=\Coind_B^{A'}\underB=\Coind_{A}^{A'}\Coind_B^A\underB=\Coind_A^{A'}S_N^A.
$
Thus we obtain an inclusion of algebras
$\End_A(S_N^A)\hookrightarrow\End_{A'}(S_{N'}^{A'}).$

\sssn
\Lemma \label{canis2}
Let $X\bul\in\dC(N)$. Then we have  a canonical isomorphism of $A'$-modules
$$
\Ind_{N'}^{A'}\circ\Coind_N^{N'}X\bul\cong\Coind_A^{A'}\circ\Ind_N^AX\bul.
$$
\dok
Suppose $X,Y\in A\mod$ and $\hom_A(X,Y)\in\dC({\cal V}ect)$. Then we have a canonical
isomorphism
$
(X\ten_AY^*)^*\cong\left(\hom_A(X,Y)\right)^*
$.
Using this isomorphism we obtain
\begin{gather*}
\Ind_{N'}^{A'}\circ\Coind_N^{N'}X\bul\til{\map}
\left(\Coind_{N'}^{A'}\left(\Coind_N^{N'}X\bul\right)^*\right)^*\til{\map}
\left(\Coind_{N'}^{A'}\Ind_N^{N'}(X\bul)^*\right)^*,\\
\Coind_A^{A'}\Ind_N^AX\bul\til{\map}
\left(\Ind_A^{A'}\left(\Ind_N^AX\bul\right)^*\right)^*\til{\map}
\left(\Ind_A^{A'}\Coind_N^A(X\bul)^*\right)^*.
\end{gather*}
Now change all the gradings to the opposite ones and use Lemma~\ref{canis1}.
\qed

\sssn
\Lemma
The constructed map provides an inclusion of algebras
$\oppA\hookrightarrow A^{\prime\sharp}$ well defined with respect to the
triangular decompositions.

\dok
The inclusion $B^{\opp}\subset\End_A(S_N^A)$ (resp.
$B^{\opp}\subset\End_{A'}(S_{N'}^{A'})$) is provided by the realisations of the
semiregular modules as coinduced ones from the regular $B$-module. Thus
the constructed map $\End_A(S_N^A)\map\End_{A'}(S_{N'}^{A'})$ identifies the
images of the algebra $B^{\opp}$.

Next by Lemma~\ref{canis2} we have the  isomorphism of functors $\dC(N)\map
\dC(A')$:
$$
\Ind_{N'}^{A'}\circ\Coind_N^{N'}\cong\Coind_A^{A'}\circ\Ind_N^A.
$$
Thus the images of two inclusions
$$
N^{\opp}\hookrightarrow\End_A(S_N^A)\hookrightarrow\End_{A'}(S_{N'}^{A'})
\text{ and }N^{\opp}\hookrightarrow N^{\prime\opp}\hookrightarrow\End_{A'}(S_{N'}^{A'})
$$
coincide. Now recall that as a vector space the subalgebra
$\oppA\subset\End_A(S_N^A)$ is described as $N^\opp\ten B^{\opp}$.\qed

Note also that the inclusions $A\subset A'$ and
$\oppA\subset A^{\prime\sharp}$ provide the same morphism of DG-algebras
$\Db(A^{\prime\sharp},B)^{\opp}=\Db(A',B)\map\Db(A,B)=\Db(\oppA,B)^{\opp}$.

\sssn  \label{mor2}
Suppose we have an $A^{\prime\sharp}$-module $M^{\sharp}$.
Consider the  morphism of functors
$$
\Res_{\oppA}^{A^{\prime\sharp}}:\
\Rhom_{A^{\prime\sharp}}(\underline{B}_{A^{\prime\sharp}},M^{\sharp})\map
\Rhom_{\oppA}(\underline{B}_{\oppA},M^{\sharp}).
$$
It is represented by a morphism
of DG-modules over  $\Db(A',B)^{\opp}$ constucted similarly to the one in
Lemma~\ref{morph}:
$\dD_{A^{\prime\sharp}}(M^{\sharp})\map\dD_{\oppA}(M^{\sharp})$.

\sssn
\Prop  \label{posext}
There exists a natural (functorial) morphism
$$
\Exts_{A}(M^{\sharp},M)\map
\Exts_{A'}(M^{\sharp},M).
$$
\dok
For any two $\Db(A,B)$-DG-modules $X\bul$ and $Y\bul$
the morphism of the DG-algebras $\Db(A',B)\map\Db(A,B)$ obtained from the morphism of functors
$\Rhom_{A'}(*,*)\map\Rhom_A(\Res_A^{A'}(*),\Res_A^{A'}(*))$ provides a natural
restriction map
$$
\Rhom_{\Db(A,B)}(X\bul,Y\bul)\map
\Rhom_{\Db(A',B)}(X\bul,Y\bul).
$$
In particular we have a natural morphism
$$
\Rhom_{\Db(A,B)}(\dD_{\oppA}(M^{\sharp}),\upD_A(M))\map
\Rhom_{\Db(A',B)}(\dD_{\oppA}(M^{\sharp}),\upD_A(M)).
$$
Now the required morphism is obtained by the composition of the constructed one with the maps
described in~\ref{mor1} and~\ref{mor2}.\qed

Similarly we obtain the following statement.
\sssn
\Prop   \label{pos}
Suppose we have an inclusion of algebras $A\subset A'$ such that both algebras have triangular decompositions
$A=B\ten N$ and $A'=B\ten N'$ as vector spaces and the inclusion preserves the triangular decompositions, i~e.
$N\subset N'$. Let $M\in\dC(A')$ and
$M^{\sharp}\in\dC(A^{\prime\sharp})$.  Then there exists a natural
(functorial) morphism $$ \Tors^{A'}(M^{\sharp},M)\map
\Tors^{A}(M^{\sharp},M).\qed
$$

\sssn
\Lemma       \label{canmap1}
let $M^{\sharp}\in\dC(A^{\prime\sharp})$. Then
the canonical map
$$
\Tors^{A'}(M^{\sharp},S_{N'}^{A'})
\map
\Tors^{A}(M^{\sharp},S_{N'}^{A'})
$$
coincides with the obvious one
$M^{\sharp}\map\Coind_{N^{\opp}}^{N^{\prime\opp}}\Res_{N^{\opp}}^{A^{\prime\sharp}}M^{\sharp}$.
\qed

\subsection{Semiinfinite induction functor.} \label{sind}
Suppose we have an inclusion of
associative algebras $A\subset A'$ such that, as before, both algebras
have triangular decompositions $A=B\ten N$ and $A'=B'\ten N'$ and the
inclusion preserves them --- $B\subset B'$  and $N\subset N'$, but the
inclusion is neither negative nor positive. Suppose also for simplicity
that $B'$ (resp. $N'$) is free both as a right and  as a left $B$-  (resp.
$N$-)module, $B'=B\ten V^-$ and $N'=N\ten V^+$; the spaces of generators
are graded with finite dimensional grading components. Then it is still
possible to construct  a morphism of $^\sharp$-algebras as follows.

\sssn
\Lemma  \label{con}
We have a canonical isomorphism of $A^{\prime\sharp}$-$A$-bimodules
$$
S_{N'}^{A'}\cong\tor_{\frac\infty 2+0}^A(S_{N'}^{A'},S_N^A).\qed
$$
Note that $\oppA$ acts on the right hand side bimodule via the second factor.
Thus we have a moprhism of algebras $\oppA\map\End_A(S_{N}^{A'})$. A simple calculation shows that the image in fact belongs to the space
of continious endomorphisms. By~\ref{contin} we obtain an inclusion of algebras
$\oppA\map A^{\prime\sharp}$. It is easily checked that the pair of algebras
$(\oppA,A^{\prime\sharp})$ satisfies the conditions~\ref{sind}.
Consider the $A'$-$A$-bimodule $\hom_{\oppA}^\c(S_{N}^{A},S_{N'}^{A'})$.

\sssn
\Lemma
$\hom_{\oppA}^\c(S_{N}^{A},S_{N'}^{A'})\til{\map} V^{+*}\ten V^-\ten A$ as a right
$A$-module.  \qed

Consider the {\em semiinfinite induction} functor
$$
\Ind^{\si}:\ A\mod\map A'\mod,\ M\mapsto
\hom_{\oppA}^\c(S_{N}^{A},S_{N'}^{A'}) \ten_AM.
$$
By the previous Lemma it
is exact and well defined on corresponding derived categories:  $\Ind^\si:
{\cal D}(A\mod)\map{\cal D}(A'\mod)$. Evidently for $M\in A\mod$ we have
$\Ind^\si(M)=V^{+*}\ten V^-\ten M$ as a vector space.

\section{$^{\sharp}$-algebras for affine quantum groups}

In this section we calculate the endomorphism algebras of semiregular modules
over affine quantum groups at a fixed level with respect to various subalgebras.
From now on we fix the base field $\k=\Q(v)$.

\ssn
Consider the affine quantum group $\U_k$ at a fixed level $k\in\Z$ (see~\ref{level}).
Recall that $\U$ is a $Y$-graded algebra. The $Y$-grading on $\U_k$ is induced by the one on $\U$.
We define a $\Z$-grading of a homogenious element $u\in\U_k$ by
$\deg u:=\hgt(\deg^Y u)$.

The triangular decomposition of $\U_k$ is provided by the one of $\U$. Namely,
in terms of the previous section, we have
$A=\U_k$, $B=\U_k^{\le0}=\U^-\ten\U_k^0$ and $N=\U^+$ where $\U_k^0$ denotes the
quotient algebra of $\U^0$ by the relation $Z=v^k$. Evidently the triple
$(\U_k,\U^{\le0}_k,\U^+)$ satisfies the conditions (i)-(v) from~\ref{setup}.

Consider the semiregular module $S_{\U^+}^{\U_k}$ over $\U_k$ with respect to $\U^+$.

\sssn
\Lemma
The algebra $\U_k$ satisfies also the condition (vi) from~\ref{setup2}.\qed

\ssn
It is known (see e. g. [V], [Ar2]) that the $^{\sharp}$-algebra for the universal enveloping algebra
of the affine Lie algebra at a fixed level $k$ is equal to the universal enveloping algebra
of the same Lie algebra at the level $-2h^{\vee}-k$ where $h^{\vee}$ denotes the dual Coxeter number
(see~\ref{fix}). Here we present a similar statement for the affine quantum groups.

\sssn
Consider the DG-algebra $\Db(\U_k,\U_k^{\le0})$. Then the multiplication in $\Db(\U_k.\U_k^{\le0})$
provides the isomorphisms of graded vector spaces
$$
\Db(\U^+,\k)\ten\U_k^{\le0}\til{\map}\Db(\U_k,\U_k^{\le0})\text{ and }
\U_k^{\le0}\ten\Db(\U^+,\k)\til{\map}\Db(\U_k,\U_k^{\le0}).
$$
Recall also that $\Db(\U^+,\k)$ is isomorphic to the tensor algebra of the graded vector space
$\overline{\U}^{+*}$. We are going to find the DG-algebra $\Db(\U,\U^{\le0})^{\opp}$ explicitly.

Consider the map $ s:\ '\til{\U}\map{}'\til{\U}^{\opp}$ defined by
$$
 s(E_i)=-E_i,\  s(F_i)=-F_i,\  s(K_\mu)=K_{-\mu},\ i\in I,\ \mu\in Y.
$$
\sssn
\Lemma
The map $s$ is well defined on $\til{\U}$ and provides
an isomorphism of algebras $\U_k\map\U_{-k}^{\opp}$.
\qed

In particular $s$ provides isomorphisms of algebras
$$
 s_{\U^{\le0}}:\ \U_k^{\le0}\til{\map}\U_{-k}^{\le0\opp}\text{ and }
 s_{\U^+}:\ \U^+\til{\map}\U^{+\opp}.
$$
Note that $s$ has nothing to do with the antipode map being a part of the Hopf algebra structure on
$\til{\U}$.

We denote the corresponding map of dual spaces by $ s_{\U^+}^*:\ \U^{+*}\map\U^{+*}$.
Consider the morphism of algebras $D( s_{\U^+}):\ T(\overline{\U}^{+*})\map T(\overline{\U}^{+*})^{\opp}$
as follows:
$$
D( s_{\U^+})(u_1\ten\ldots\ten u_m):= s_{\U^+}^*(u_m)\ten\ldots\ten s_{\U^+}^*(u_1),\
u_1,\ldots,u_m\in \overline{\U}^+.
$$
\sssn
\Lemma          \label{positive}
The map $D( s_{\U^+})$ preserves the differentials and defines an isomorphism
of DG-algebras $\Db(\U^+,\k)\til{\map}\Db(\U^+,\k)^{\opp}$.
\qed

Now we use the antipode map in $\U$ and $\U_k$ to perform a similar construction
for these algebras.

\subsubsection{Antipode maps in DG-algebras.} We prefer to work in a more general setup.
Suppose we have an associative algebra $A$ with two its subalgebras $B$ and $N$
like in~\ref{setup} with a triangular decomposition satisfying the
conditions (i)-(vi).  Suppose also there exists an antiautomorphism $s_A:\
A\til{\map}A^{\opp}$ preserving the triangular decomposition, i.~e. $
s_A(B)\subset B$ and $ s_A(N)\subset N$, for simplicity we set
$s_A^2=\operatorname{Id}$.  By~\ref{free} the algebra $\Db(A,B)$ is
isomorphic to the algebra $T_B(M)$ for the $B$-bimodule
$M=\hom_{B\text{-left}}(A/B,B)$.

\sssn \label{homhom}
\Lemma
The map $D(s_A):\ \hom_{B\text{-left}}(A/B,B)\map\hom_{B\text{-right}}
(A/B,B)$ given by
$$
D(s_A)(f)(a)=s_B(f(s_A(a))),\ f:\ A/B\map B,\ a\in A,
$$
provides a correctly defined isomorphism of $B$-bimodules. Here the left
$B$-module structure on the space $\hom_{B\text{-right}}(A/B,B)$ is
provided by the {\em left} $B$-multiplication in $A/B$ twisted by $s_A$
and the right $B$-module structure on it is provided by the {\em
left} $B$-multiplication in $B$ also twisted by $s_B$.

\dok
First we check that $D(s_A)(f)$ belongs to the described space of homomorphisms:
\begin{gather*}
D(s_A)(f)(ab)=s_B(f(s_A(ab)))=s_B(f(s_B(b)s_A(a)))\\=
s_B(s_B(b)f(s_A(a)))=s_B(f(s_A(a)))s_B^2(b)=D(s_A)(f)(a)b.
\end{gather*}
Next we check that $D(s_A)$ is a morphism of the left $B$-modules:
\begin{gather*}
D(s_A)(b\cdot f)(a)=s_B(b\cdot f(s_A(a)))=s_B(f(s_A(a)b))\\=
s_B(f(s_A(s_B(b)a)))=D(s_A)(f)(s_B(b)a)=(b\cdot D(s_A)(f))(a).
\end{gather*}
Finally we check that $D(s_A)$ is a morphism of right $B$-modules:
\begin{multline*}
D(s_A)(f\cdot b)(a)=s_B(f\cdot b(s_A(a)))=s_B(f(s_A(a))b)=
s_B(b)s_B(f(s_A(a)))\\=(D(s_A)(f)\cdot b)(a).\qed
\end{multline*}
Consider the $B^{\opp}$-bimodule $M^{\opp}:=\hom_{B^{\opp}\text{-left}}
(A^{\opp}/B^{\opp},B^{\opp})$ with the left $B^{\opp}$-module structure
provided by the {\em right} $B^{\opp}$-multiplication in $A^{\opp}/B^{\opp}$
(i.~e. by the {\em left} $B$-multiplication in $A/B$) and the right
$B^{\opp}$-module structure provided by the {\em right} $B^{\opp}$-multiplication
in $B^{\opp}$ (i.~e. by the {\em left} $B$-multiplication in $B$). Consider as before
the $B^{\opp}$-free algebra $T_{B^{\opp}}(M^{\opp})$ generated by $M^{\opp}$.

We extend the moprhism $D(s_A)$ to the following map
\begin{gather*}
D(s_A):\ T_B(M)\map T_{B^{\opp}}(M^{\opp}),\
b\mapsto s_B(b);\\
\text{(*)}\qquad \qquad f_1\ten\ldots\ten f_m\mapsto D(s_A)(f_m)\ten\ldots\ten
D(s_A)(f_1).
\end{gather*}
\sssn
\Lemma
The map $D(s_A)$ provides a correctly defined isomorphism of associative algebras $T_B(M)^{\opp}$
and $T_{B^{\opp}}(M^{\opp})$.

\dok
Note that $\hom_{B\text{-right}}(A/B,B)$ is isomorphic to the $B$-bimodule obtained
from $M^{\opp}$ by composing the $B$-actions with the antipode map $s_B$.
Thus by Lemma~\ref{homhom} we have a vector space map $M\map M^{\opp}$.
Any homomorphism of vector spaces $V_1\map V_2$ provides an antihomomorphism of tensor algebras
(over the base field) $T(V_1)\map T(V_2)$ given by the formula (*). We are to check that
the map is well defined with respect to the $B$-action. The nessesary calculation looks as follows:
\begin{gather*}
D(s_A)((f_1\cdot b)\ten f_2)=
D(s_A)(f_2)\ten D(s_A)((f_1\cdot b)=
D(s_A)(f_2)\ten D(s_A)(f_1)\cdot s_B(b)\\=
D(s_A)(f_2)\ten s_B(b)\cdot^{\opp}D(s_A)(f_1)=
D(s_A)(f_2)\cdot^{\opp}s_B(b)\ten D(s_A)(f_1)\\=
s_B(b)\cdot D(s_A)(f_2)\ten D(s_A)(f_1)=
D(s_A)(b\cdot f_2)\ten D(s_A)(f_1)=
D(s_A)(f_1\ten b\cdot f_2).
\end{gather*}
Here $\cdot^{\opp}$ denotes the left and right actions of $B^{\opp}$ on $M^{\opp}$
to distinguish from $\cdot$ denoting the same actions considered as the right and left actions
of $B$ respectively.
\qed

Now note that by~\ref{free} the associative algebra $T_{B^{\opp}}(M^{\opp})$
is iosmorphic to the algebra $\Db(A^{\opp},B^{\opp})$. Thus we have the isomorphism of algebras
(with differentials forgotten)
$D(s_A):\ \Db(A,B)^{\opp}\til{\map}\Db(A^{\opp},B^{\opp})$.

\sssn
\Cor
We have the isomorphisms of associative algebras

\qquad \text{(i)}$\Db(\U,\U^{\le0})^{\opp}\til{\map}\Db(\U,\U^{\le0})$;

\qquad \text{(ii)}$\Db(\U_k,\U_k^{\le0})^{\opp}\til{\map}\Db(\U_{-k},\U_{-k}^{\le0})$.\qed

Again, these isomoprhisms do not nessesarily preserve the differentials in the DG-algebras.
Note that the associative algebra structure in $\Db(\U_k,\U_k^{\le0})$ does not depend on the
level $k$, and it is only the differential that ``remembers'' the level.

\sssn
\Prop
The map $D(s_{\U})$ provides an isomorphism of the DG-algebras
$\Db(\U_k,\U_k^{\le0})^{\opp}\til{\map}\Db(\U_{-2h^{\vee}-k},\U_{-2h^{\vee}-k}^{\le0})$.

\dok
Since the associative algebras $\Db(\U_k,\U_k^{\le0})$ do not depend on $k$, we
can identify them and think of the whole picture as of a family of differentials
$\{d_k\}$ on the graded algebra, say, $\Db(\U_0,\U_0^{\le0})$. We are to check that
the differential defined by $D(s_{\U_0})\circ d_k\circ D(s_{\U_0})^{-1}$ coincides with
$d_{-2h^{\vee}-k}$. We fix the isomorphism of vector spaces
$$
\Db(\U_0,\U_0^{\le0})\til{\map}\U_0^{\le0}\ten\Db(\U^+,\k)=\U_0^{\le0}\ten T(\overline
{\U}^{+*}).
$$
Recall that $T(\overline{\U}^{+*})$
with the differential generated by the map dual to the multiplication map in
$\U^+$ is a DG-subalgebra in all the DG-algebras $\Db(\U_k,\U_k^{\le0})$.
By Lemma~\ref{positive} the restriction of
the differential $D(s_{\U_0})\circ d_k\circ D(s_{\U_0})^{-1}$ to $T(\overline{\U}^{+*})$
coincides with the discribed one. On the other hand the associative algebra
$\Db(\U_0,\U_0^{\le0})$ is generated by $\U_0^{\le0}$ and $T(\overline{\U}^{+*})$
thus it is sufficient to check that the components $\U_0^{\le0}\map\U_0^{\le0}\ten\overline{\U}^{+*}$
of $d_{-2h^{\vee}-k}$ and $D(s_{\U_0})\circ d_k\circ D(s_{\U_0})^{-1}$ coincide.
In fact, using the Leibnitz rule, we can restrict  the differentials to the space of generators
of $\U_0^{\le0}$ or, since $\U_0^{\le0}\ten\overline{\U}^{+*}$ is strictly negatively graded and the differentials
preserve the gradings, we have to check only that
$$
d_{-2h^{\vee}-k}(F_i)=D(s_{\U_0})\circ d_k\circ D(s_{\U_0})^{-1}(F_i)\text{ for }i\in I.
$$
We omit the direct calculation that looks completely like the one in [Ar2],
Corollary 4.4.2, in the case of affine Lie algebras.
\qed

\sssn
\Cor
There exists an inclusion of associative algebras $\U_{-2h^{\vee}-k}\subset
\End_{\U_k}(S_{\U^+}^{\U_k})$.\qed

In particular $S_{\U_k}$ becomes a $\U_k$-$\U_{-2h^{\vee}-k}$-bimodule.

\ssn
Fix a positive weight $x\in Y^{\prime\prime+}$, $\langle i,x\rangle>0$ for every $i\in \overline{I}$,
and the  transvection element $\theta_{mx}\in T\subset W$.
Recall that we have defined the subalgebra $\U_{\theta_{mx}}^+\subset\U^+$.
Note  that repeating word by word the considerations from the previous subsection
we obtain a proof of the following statement.

\sssn
\Cor   \label{theta}
There exists an inclusion of associative algebras $\U_{-2h^{\vee}-k}\subset
\End_{\U_k}(S_{T_{\theta_{mx}}^{-1}(\U^-)}^{\U_k})$.\qed

Consider the semiregular $\U_k$-module $S_{\U_{\theta_{mx}}^+}^{\U_k}:=\Coind_{\U_{\theta_{mx}}^+}^{\U_k}\Ind
_{\k}^{\U_{\theta_{mx}}^+}\underk$
with respect to the subalgebra $\U_{\theta_{mx}}^+$. Here we prove that like in the affine Lie algebra case
there exists an associative algebra inclusion $\U_k\hookrightarrow\End_{\U_k}
(S_{\U_{\theta_{mx}}^+}^
{\U_k})$. Unfortunitely unlike in the Lie algebra case we cannot write down the second action
of $\U_k$ on the semiregular module explicitly.
Instead of that we prove the following statement.

\sssn
\Theorem  \label{bimodule}
There exists a natural isomorphism in the derived category of left $\U_k$-modules
$$
S_{\U^+}^{\U_k}
\tenl_{\U_{-2h^{\vee}-k}}
S_{T_{\theta_{mx}}^{-1}(\U^-)}^{\U_{-2h^{\vee}-k}}
\til{\map}S_{\U_{\theta_{mx}}^+}^{\U_k}[\lth(\theta_{mx})].
$$
Here $[\cdot]$ denotes the shift in the derived category and the left $\U_k$-module
structure on the left hand side of the isomorphism is provided by the natural left
$\U_k$-module structure on $S_{\U^+}^{\U_k}$.

\Rem
In other words we have
$$
\tor_{\ne\lth(\theta_{mx})}^{\U_k}
(S_{\U^+}^{\U_k},S_{T_{\theta_{mx}}^{-1}(\U^+)}^{\U_{-2h^{\vee}-k}})=0\text{ and }
\tor_{\lth(\theta_{mx})}^{\U_k}
(S_{\U^+}^{\U_k},S_{T_{\theta_{mx}}^{-1}(\U^+)}^{\U_{-2h^{\vee}-k}})=
S_{\U_{\theta_{mx}}^+}^{\U_k}
$$
as a left $\U_k$-module. Now note that by Corollary~\ref{theta} the left hand side of the latter equivalence
carries the natural structure of a $\U_k$-bimodule.

\subsection{Proof of Theorem~\ref{bimodule}.}
The statement of the Theorem will be derived from Lemmas~\ref{a}--\ref{e}.
Consider the natural triangular decomposition of the algebra $T_{\theta_{mx}}^{-1}(\U^-)$
provided by its inclusion into $\U$.

\sssn \label{a}
\Lemma
The algebra $T_{\theta_{mx}}^{-1}(\U^-)^{\sharp}$ is isomorphic to
$T_{\theta_{mx}}^{-1}(\U^-)$.

\dok
By~\ref{sind} the algebra $T_{\theta_{mx}}^{-1}(\U^-)$ is a subalgebra in
$\U_k^{\sharp}$ generated by $\U_{\theta_{mx}}^+$ and $T_{\theta_{mx}}^{-1}(\U^-)\cap
\U^-$. On the other hand we know by Theorem~\ref{bimodule} that
$\U_k^{\sharp}\cong\U_{-2h^{\vee}-k}$.  Finally, the subalgebra in
$\U_{-2h^{\vee}-k}$ generated by $\U^+_{\theta_{mx}}$ and
$T_{\theta_{mx}}^{-1}(\U^-)\cap\U^-$ is evidently isomorphic to
$T_{\theta_{mx}}^{-1}(\U^-)$.  \qed

The following statement is a particular case of Lemma~\ref{con}.

\sssn
\Lemma
There exists a canonical isomorphism of $\U_k$-$T_{\theta_{mx}}^{-1}(\U^-)$-bimodules
$$
S_{\U^+}^{\U_k}\til{\map}
\Ind_{T_{\theta_{mx}}^{-1}(\U^-)}^{\si\quad\U_k}
S_{\U_{\theta_{mx}}^+}^{T_{\theta{mx}}^{-1}(\U^-)}.\qed
$$
\sssn
\Lemma
\begin{multline*}
\left(\Ind_{T_{\theta_{mx}}^{-1}(\U^-)}^{\si\quad\U_k}
S_{\U_{\theta_{mx}}^+}^{T_{\theta{mx}}^{-1}(\U^-)}\right)\tenl_{
T_{\theta_{mx}}^{-1}(\U^-)}
T_{\theta_{mx}}^{-1}(\U^-)^*
\til{\map}\\
\Ind_{T_{\theta_{mx}}^{-1}(\U^-)}^{\si\quad\U_k}
\left(S_{\U_{\theta_{mx}}^+}^{T_{\theta{mx}}^{-1}(\U^-)}\tenl_{
T_{\theta_{mx}}^{-1}(\U^-)}
T_{\theta_{mx}}^{-1}(\U^-)^*\right)
.
\end{multline*}
\dok
The statement follows from the exactness of the semiinfinite induction functor.
\qed

\sssn
Recall that a nonnegatively graded algebra $A=\underset{n\ge0}{\bigoplus}A_n$,
$A_0=\k$, is called {\em quadratic} if it is generated by the space $A_1$, $\dim A_1<\infty$,
and its relations ideal in $T(A_1)$ is generated by a space $J\subset A_1\ten A_1$.
The {\em quadratic dual algebra} $A^!$ for the quadratic algebra $A$ is by definition
the algebra on generators $A_1^*$ and quadratic relations $J^{\perp}\subset A_1^*\ten A_1^*$.
The definitions of the Koszul  complex  $K_A\bul$ (resp. the co-Koszul complex
$C_A\bul$) for the algebra $A$ can be found e.~g. in [Pr]. Note only that
$K_A\bul=A\ten A^{!*}$ and $C_A\bul=A\ten A^!$ as vector spaces. A quadratic algebra
$A$ is called {\em Koszul} (resp. {\em co-Koszul}) if its Koszul (resp. co-Koszul)
complex is quasiisomorphic to the trivial $A$-module $\underk$.

Note that the algebra $A^!$ for the Koszul algebra $A$ is equal to
$\Ext_A(\underk,\underk)$ since the Koszul complex provides a $A$-free
resolution of the trivial module.  Evidently $A$ is Koszul iff $A^!$ is
so.

By definition the {\em twisted exterior algebra} $\Lambda_q$ is given by the generators
$E_\beta^*,\ \beta\in R_{\theta_{mx}}^+$ and relations
\begin{gather*}
E_\alpha^*E_\beta^*+E_{\beta}^*E_\alpha^*=0\text{ if }\beta>\alpha;\\
E_\alpha^{*2}=0\text{ for }\alpha\in R_{\theta_{mx}}^+.
\end{gather*}
\sssn
\Lemma  \label{b}
$\gr\U_{\theta_{mx}}^+$ is a quadratic Koszul algebra with the quadratic dual algebra
equal to $\Lambda_q$.\qed

Recall that a finite dimensional algebra $A$ is called {\em strongly Frobenius}
if the coregular $A$-bimodule $A^*$ is isomorphic to the regular $A$-bimodule.

\sssn
\Lemma           \label{b1}
Suppose that  a finite dimensional quadratic Koszul algebra $A$ is strongly Frobenius.
Then it is also co-Koszul.

\dok
The isomorphism $A^*\cong A$ identifies the complexes $C_A\bul$ and
$K_A^{\bullet*}$.  \qed

Note that the algebra $\Lambda_q$ is strongly Frobenius with the required bimodules'
isomorphism given by the pairing $\Lambda_q\times \Lambda_q\map\k$
that is defined as a projection of the product in the algebra on
the top grading component.

\sssn
\Cor
$\Lambda_q$ is both Koszul and co-Koszul.

\dok
follows immediately from Lemmas~\ref{b}~and~\ref{b1}.\qed

\sssn
\Lemma
$\tor_{\ne\lth(\theta_{mx})}^{\gr\U_{\theta_{mx}}^+}
(\underk,\U_{\theta_{mx}}^{+*})=0,\
\tor_{\lth(\theta_{mx})}^{\gr\U_{\theta_{mx}}^+}
(\underk,\U_{\theta_{mx}}^{+*})=\k.$

\dok
We calculate the $\tor$ spaces using the Koszul resolution of $\underk$.
We have
$$
K_{\gr\U_{\theta_{mx}}^+}\bul\ten_
{\gr\U_{\theta_{mx}}^+}
\gr\U_{\theta_{mx}}^{+*}=
\Lambda_q^*\ten
\gr\U_{\theta_{mx}}^{+*}=
C_{\Lambda_q}^{\bullet*}[\lth(\theta_{mx})].
$$
Now use the previous Lemma.
\qed

A similar satatement holds evidently for $\gr\U_{\theta_{mx}}^-$.

\sssn
\Lemma  \label{c}

\qquad\text{(i)}
$\Tor^{\U_{\theta_{mx}}^{\pm}}(\underk,\U_{\theta_{mx}}^{\pm*})
=\delta_{\bullet,\lth(\theta_{mx})}\k$;\\

\qquad\text{(ii)} for any filtered graded $\U_{\theta_{mx}}^{\pm}$-module
$M=\underset{n\ge0}{\bigcup}F^nM,\ \dim F^nM/F^{n-1}M<\infty$,
$$
\Tor^{\U_{\theta_{mx}}^{\pm}}(M,\U_{\theta_{mx}}^{\pm*})
=\delta_{\bullet,\lth(\theta_{mx})}M;\\
$$
\qquad\text{(iii)} $\Tor^{\U_{\theta_{mx}}^{\pm}}(\U_{\theta_{mx}}^{\pm*},
\U_{\theta_{mx}}^{\pm*})
=\delta_{\bullet,\lth(\theta_{mx})}\U_{\theta_{mx}}^{\pm*}$.

\dok
We prove the  statements of the Lemma for the algebra $\U_{\theta_{mx}}^+$.
The case of the algebra $\U_{\theta_{mx}}^-$  is treated similarly.

The well known spactral sequence
$$
\Tor^{
\gr\U_{\theta_{mx}}^+(\bullet)}
(\underk,
\gr\U_{\theta_{mx}}^{+*})
)\Longrightarrow
\Tor^{
\U_{\theta_{mx}}^+}
(\underk,
\U_{\theta_{mx}}^{+*})
$$
shows that
$\underk\tenl_
{\U_{\theta_{mx}}^+}
\U_{\theta_{mx}}^{+*}$
itself is isomorphic to
$\U_{\theta_{mx}}^{+*}[\lth(\theta_{mx})]$.
Here $(\bullet)$ denotes the second grading on $\tor$ spaces obtained
from the $\Z_{\ge0}^{R_{\theta_{mx}}^+}$-grading. Thus the first statement of the Lemma is proved.

To prove the second one consider the standard complex for the computation of
$\tor$ spaces as follows:
$$
M\ten_{\U_{\theta_{mx}}^+}\Barb(\U_{\theta_{mx}}^+,\k,\U_{\theta_{mx}}^{+*})=
M\ten T(\overline{\U}_{\theta_{mx}}^+)\ten\U_{\theta_{mx}}^{+*}.
$$
Note that $M=\dirlim F^nM$ and $
M\ten T(\overline{\U}_{\theta_{mx}}^+)\ten\U_{\theta_{mx}}^{+*}=\dirlim
F^nM\ten T(\overline{\U}_{\theta_{mx}}^+)\ten\U_{\theta_{mx}}^{+*}$.
Since $\dirlim$ functor is exact we have
$$
\Tor^{\U_{\theta_{mx}}^{\pm}}(M,\U_{\theta_{mx}}^{\pm*})=
\dirlim\Tor^{\U_{\theta_{mx}}^{\pm}}(F^nM,\U_{\theta_{mx}}^{\pm*})
=\delta_{\bullet,\lth(\theta_{mx})}\dirlim F^nM=
=\delta_{\bullet,\lth(\theta_{mx})}M.
$$
The third statement is an immidiete concequence of the second one since
the filtration $F^n\U_{\theta_{mx}}^{+*}=\underset{k\ge -n}{\bigoplus}
\left(\U_{\theta_{mx}}^{+*}\right)_k$ satisfies the condition in (ii).
\qed

\sssn
\Lemma \label{d}
$
S_{\U_{\theta_{mx}}^+}^{T_{\theta{mx}}^{-1}(\U^-)}\tenl_{
T_{\theta_{mx}}^{-1}(\U^-)}
T_{\theta_{mx}}^{-1}(\U^-)^*
\til{\map}
T_{\theta_{mx}}^{-1}(\U^-)^*[\lth(\theta_{mx})]$
in the derived category of left graded
$T_{\theta_{mx}}^{-1}(\U^-)$-modules.

\dok
First we construct the required morphism
in the derived category of left graded
$T_{\theta_{mx}}^{-1}(\U^-)$-modules.

By Shapiro Lemma we have
\begin{gather*}
\Rhom_{\U_{\theta_{mx}}^+}(
S_{\U_{\theta_{mx}}^+}^{T_{\theta_{mx}}^{-1}(\U^-)}\tenl_{
T_{\theta_{mx}}^{-1}(\U^-)}
T_{\theta_{mx}}^{-1}(\U^-)^*,
T_{\theta_{mx}}^{-1}(\U^-)^*)\\=
\Rhom_{\k}
(S_{\U_{\theta_{mx}}^+}^{T_{\theta_{mx}}^{-1}(\U^-)}\tenl_{
T_{\theta_{mx}}^{-1}(\U^-)}
T_{\theta_{mx}}^{-1}(\U^-)^*,
\k)\\=
\left(
S_{\U_{\theta_{mx}}^+}^{T_{\theta_{mx}}^{-1}(\U^-)}\tenl_{
T_{\theta_{mx}}^{-1}(\U^-)}
T_{\theta_{mx}}^{-1}(\U^-)^*\right)^*
=\left(
\U_{\theta_{mx}}^{+*}\tenl_
{\U_{\theta_{mx}}^+}
T_{\theta_{mx}}^{-1}(\U^-)^*\right)^*.
\end{gather*}
By~\ref{dec} we have
$
T_{\theta_{mx}}^{-1}(\U^-)=
\U_{\theta_{mx}}^+\ten T_{\theta_{mx}}^{-1}(\U^-)\cap\U^-
$.
Thus by the previous Lemma we obtain
\begin{gather*}
\Rhom_{\U_{\theta_{mx}}^+}(
S_{\U_{\theta_{mx}}^+}^{T_{\theta{mx}}^{-1}(\U^-)}\tenl_{
T_{\theta_{mx}}^{-1}(\U^-)}
T_{\theta_{mx}}^{-1}(\U^-)^*,
T_{\theta_{mx}}^{-1}(\U^-)^*)\\=
\left(
\U_{\theta_{mx}}^{+*}\tenl_
{\U_{\theta_{mx}}^+}
\U_{\theta_{mx}}^{+*}\right)^*
\ten
T_{\theta_{mx}}^{-1}(\U^-)\cap\U^-\\=
\U_{\theta_{mx}}^+\ten
T_{\theta_{mx}}^{-1}(\U^-)\cap\U^-
[-\lth(\theta_{mx})]
=T_{\theta_{mx}}^{-1}(\U^-)
[-\lth(\theta_{mx})].
\end{gather*}
In particular the canonical element in
$$
\Rhom_{\U_{\theta_{mx}}^+}(
S_{\U_{\theta_{mx}}^+}^{T_{\theta{mx}}^{-1}(\U^-)}\tenl_{
T_{\theta_{mx}}^{-1}(\U^-)}
T_{\theta_{mx}}^{-1}(\U^-)^*,
T_{\theta_{mx}}^{-1}(\U^-)^*)
$$
that corresponds to $1\in
T_{\theta_{mx}}^{-1}(\U^-)
[-\lth(\theta_{mx})]$
provides the required morphism in the derived category of left graded
$T_{\theta_{mx}}^{-1}(\U^-)$-modules.

Finally we prove that this element is an
isomorphism in the derived category.  Since the morphism in the derived
category is already constructed, we need only to check that it is a
quasiisomorphism on the level of complexes of vector spaces. The
calculation almost repeats the previous one:
\begin{gather*}
S_{\U_{\theta_{mx}}^+}^{T_{\theta_{mx}}^{-1}(\U^-)}\tenl_
{T_{\theta_{mx}}^{-1}(\U^-)}T_{\theta_{mx}}^{-1}(\U^-)^*
\til{\map} \U_{\theta_{mx}}^{+*} \tenl
_{\U_{\theta_{mx}}^+} T_{\theta_{mx}}^{-1}(\U^-)^*\\ \til{\map}
\U_{\theta_{mx}}^{+*} \tenl _{\U_{\theta_{mx}}^+}
\left(\U_{\theta_{mx}}^{+*} \ten (\U^-\cap
T_{\theta_{mx}}^{-1}(\U^-))^*\right)
\\
\til{\map}
\left(\U_{\theta_{mx}}^{+*} \tenl _{\U_{\theta_{mx}}^+}
\U_{\theta_{mx}}^{+*}\right) \ten (\U^-\cap T_{\theta_{mx}}^{-1}(\U^-))^*
\til{\map}
T_{\theta_{mx}}^{-1}(\U^-)^* [\lth(\theta_{mx})].\qed
\end{gather*}
\sssn
\Lemma  \label{e}
$
\Ind_{T_{\theta_{mx}}^{-1}(\U^-)}^{\si\quad\U_k}
T_{\theta_{mx}}^{-1}(\U^-)^*
\til{\map}
\Coind_{\U_{\theta_{mx}}^+}^{\U_k}\U_{\theta_{mx}}^+=S_{\U_{\theta_{mx}}^+}^{\U_k}
$.\qed

By~\ref{sind} the semiinfinite induction functor is well defined on corresponding derived
categories, thus it takes the quasiisomorphism from Lemma~\ref{d} to a quasiisomorphism.
Theorem~\ref{bimodule} is proved.\qed

\section{Quantum twisted Verma modules and semiinfinite homology of
the algebra $T_{\theta_{mx}}(\U^-)$.}

In this section we show that the semiinfinite homology space
of the trivial module over $T_{\theta_{mx}}(\U^-)$
has a base enumerated by elements of the affine Weyl group
graded by the  $\theta_{mx}^{-1}$-twisted length function on $W$ (see~\ref{tlength})
just like in the affine Lie algebra case (see [Ar3]).

\subsection{Categories of $\U$-modules.}
Consider the category $\M$
of $X\times\Z$-graded $\U$-modules $M=\underset{\lambda\in X,t\in\Z}
{\bigoplus}M_{\lambda,t}$ such that

(i) for every $i\in I$ the   standard  generators
$
E_i:\ M_{\lambda,t}\map M_{\lambda+i',t+1},\
F_i:\ M_{\lambda,t}\map M_{\lambda-i',t-1};
$\\
(ii) every $K_{\mu}\in \U^0$ acts on $M_{\lambda,t}$ by scalar
$v^{\langle \mu,\lambda\rangle}$.

Morphisms in $\M$ are  morphisms of $\U$-modules that preserve
$X\times\Z$-gradings.

Fix a nonnegative integer $k\in\Z_{\ge 0}$.
Consider a full subcategory $\M_k$ in $\M$ of modules $M$ such that $Z$ acts
by the scalar $v^k$ on $M$.

\sssn
We define the character of a $\U$-module $M\in\M$ such that $\dim
M_{\lambda,t}<\infty$ for all $\lambda\in X,t\in\Z$ by
$
\ch M:=\sum_{\lambda\in X,t\in\Z}\dim M_{\lambda,t}e^{\lambda}q^t.
$
Here $q$ is a formal variable and $e^{\lambda}$ is a formal expression.

\sssn
As usual let $\O$ denote the category of $X$-graded
$\U$-modules $M=\underset{\lambda\in X}{\bigoplus}M_{\lambda}$ such that

(i) for every $i\in I$ the   standard generators
$
E_i:\ M_{\lambda}\map M_{\lambda+i'},\
F_i:\ M_{\lambda}\map M_{\lambda-i'};
$\\
(ii) every $K_{\mu}\in\U^0$ acts on $M_{\lambda}$ by scalar
$v^{\langle \mu,\lambda\rangle}$;\\
(iii) $\dim M_{\lambda}<\infty$ for all $\lambda\in X$;\\
(iv) there exist $\lambda_1,\ldots,\lambda_n\in\h^*$ such that
$$
M_{\mu}\ne0\text{ only when }\mu\in\lambda_1+R^-\cup\ldots\cup\lambda_n+R^-.
$$
\subsection{Affine BGG resolution.}
Recall that the subalgebra $\U^+$ is normal in $\U^{\ge0}$,
that is, the left and the right ideals in $\U^{\ge0}$ generated by $\overline{\U}^+$
coincide,  and the quotient algerba of $\U^{\ge0}$ by any of these ideals
denoted by $\U^{\ge0}//\U^+$ equals to $\U^0$. For the one
dimensional $\U^0$-module $\k(\lambda)$, $\lambda\in X$,
generated by $v_\lambda$ (such that
$K_\mu\cdot v_\lambda=v^{\langle\mu,\lambda\rangle}v_{\lambda}$)
the  {\em (quantum) Verma module} $\U_k\ten_{\U^{\ge0}}
\k(\lambda)$ is denoted by $M(\lambda)$.
Evidently $M(\lambda)$ is free as a $\U^-$-module and
belongs to $\cal O$ and to $\M_k$ if we put the
$\Z$-grading of the highest weight
vector $v_{\lambda}\in M(\lambda)$ equal to zero.

For a fixed positive integral level $k$  choose a dominant weight
$\lambda\in X_k^+$.
Consider the unique simple quotient module of $M(\lambda)$ denoted by
$L(\lambda)$. Then there
exists a left resolution of $L(\lambda)$ that consists of direct sums of  Verma
modules as follows.

\sssn
\Theorem
(see e~.g. [M], Theorem 3.3)
There exists a resolution
$B\bul(\lambda)$
of $L(\lambda)$ in $\O$ and in $\M_k$  of the form
\begin{multline*}
\ldots
\map
\underset{{w\in W,}\atop{\lth(w)=m}}{\bigoplus}
M(w\cdot\lambda)\langle-\hgt(\lambda-w\cdot\lambda)\rangle
\map\ldots\\ \map
\underset{{w\in W,}\atop{\lth(w)=1}}{\bigoplus}
M(w\cdot\lambda)\langle-\hgt(\lambda-w\cdot\lambda)\rangle
\map
M(\lambda)\map L(\lambda)\map 0.
\end{multline*}
Here as usual $\langle\cdot\rangle$ denotes the shift of the $\Z$-grading
in $\M$.
\qed

Like in the affine Lie algebra case  it is known that the component of the differential
$d_{w_1,w_2}:\ B^{-k}(\lambda)\map B^{-k+1}(\lambda),\ M(w_1\cdot\lambda)\map
M(w_2\cdot\lambda)$, is nonzero iff $w_1\ge w_2$ in the Bruhat order on $W$
(see [M], Theorem 3.2).
\subsection{Twisting functors.}
We introduce the functors of the twist of the $X$-gradings.
As before, fix a positive weight $x\in Y^+$ and the  transvection element $\theta_{mx}\in T\subset W$.
Recall that we have defined the associative algebra automorphism
$$
T_{\theta_{mx}}:\ \U\til{\map}\U,\ \U_{\alpha}\map\U_{\theta_{mx}(\alpha)}.
$$
We define $T_{\theta_{mx}}:\ \M_k\map\M_k$ as follows.
For $M\in\M_k$ set
\begin{gather*}
T_{\theta_{mx}}(M)=\underset{\lambda\in
X,t\in\Z}{\bigoplus}T_{\theta_{mx}} (M)_{\lambda,t},\
T_{\theta_{mx}}(M)_{\lambda,t}:=M_{\theta_{mx}(\lambda),t+\hgt(\theta_{mx}
(\lambda)-\lambda)},\\
e_{\alpha}\in\U_{\alpha},\ v_{\lambda}\in T_{\theta_{mx}}(M)_{\lambda,t}\text{ then }
e_{\alpha}\cdot v_{\lambda}:=T_{\theta_{mx}}(e_{\alpha})(v_{\lambda})\in
M_{\theta_{mx}(\lambda+\alpha),t+\hgt (\theta_{mx}(\lambda+\alpha)-\lambda)}
\\=M_{\theta_{mx}(\lambda+\alpha),t+\hgt
(\theta_{mx}(\lambda+\alpha)-\lambda-\alpha)+\hgt\alpha}
=T_{\theta_{mx}}(M)_{\lambda+\alpha,t+\hgt\alpha}.
\end{gather*}
Evidently $T_{\theta_{mx}}$ is an equivalence of categories.

\sssn
Consider the semiregular $\U_k$-module
$
S_{\U_{\theta_{mx}}^-}^{\U_k}=\Coind_{\U_{\theta_{mx}}^-}^{\U_k}\U_{\theta_{mx}}^-
$
with respect to $\U_{\theta_{mx}}^-$.
Then by \ref{bimodule} there exists an inclusion of algebras $
\U_k \hookrightarrow \End_{\U_k}(S_{\U_{\theta_{mx}}^-}^{\U_k})$.
Consider a functor
$$
S_{mx}:\ \M_k\map \U_k\mod,\ S_{mx}(M):=
\left(\hom_{\U_k}^\c(M,
S_{\U_{\theta_{mx}}^-}^{\U_k})\right)^*.
$$
The left $\U_k$-module structure on $S_{mx}(M)$ is provided by the right $\U_k$-multiplication
in $S_{\U_{\theta_{mx}}^-}^{\U_k}$.

\sssn
\Lemma
$S_{mx}$ defines a functor $\M_k\map\M_k$.
\qed

\sssn
\Def
For $x\in Y^{\prime\prime+}$ the functor of twist by $\theta_{mx}\in T\subset W$
$$
\Phi_{mx}:=T_{\theta_{mx}}\circ S_{mx},\ \Phi_{mx}:\ \M_k\map\M_k.
$$
Let us describe the image of a  Verma module under $\Phi_{mx}$.

\sssn
\Lemma            \label{weight}
$\ch \Phi_{mx}(M(\lambda))=\ch M(\theta_{mx}\cdot\lambda)\langle
-\hgt(\lambda-\theta_{mx}\cdot\lambda)\rangle$.
\qed

In particular the highest weight vector $1\ten
v_{\lambda}\in\Phi_{mx}(M(\lambda))$ has the weight $\theta_{mx}\cdot\lambda$.

\sssn
\Def
We call the $\U_k$-module
$M_{\theta_{mx}}(\theta_{mx}\cdot\lambda):=\Phi_{mx}(M(\lambda))\langle
\hgt(\lambda-\theta_{mx}\cdot\lambda)\rangle$ the {\em quantum twisted Verma module} of the
weight $\theta_{mx}\cdot\lambda$.

\sssn
\Lemma    \label{twistedsem}
When restricted to
$T_{\theta_{mx}}(\U^-)$
the quantum twisted Verma module
$M_{\theta_{mx}}(\lambda)$ is isomorphic to the semiregular
$T_{\theta_{mx}}(\U^-)$-module
$S_
{T_{\theta_{mx}}(\U_{\theta_{mx}}^-)}
^{T_{\theta_{mx}}(\U^-)}
$.

\dok
We prove first that the right
$\U^-$-modules
$\operatorname{Res}_{\U^-}^{\U_k}
\left(\hom_{\U_k}^\c(M(\lambda),
S_{\U_{\theta_{mx}}^-}^{\U_k})\right)$
and $\Coind
_{\U_{\theta_{mx}}^-}^{\U^-}
\U_{\theta_{mx}}^-$ are isomorphic to each other.
On the level of vector spaces the calculation looks as follows.
\begin{gather*}
\text{(*)}\qquad\hom_{\U_k}^\c(M(\lambda),
S_{\U_{\theta_{mx}}^-}^{\U_k})=
\hom_{\U_k}^\c(M(\lambda),
\Coind_{\U^-}^{\U^k}\Coind_{\U_{\theta_{mx}}^-}^{\U^-}\U_{\theta_{mx}}^-)\\ =
\hom_{\U^-}^\c(M(\lambda),
\Coind_{\U_{\theta_{mx}}^-}^{\U^-}\U_{\theta_{mx}}^-)=
\hom_{\U^-}^\c(\U^-,
\Coind_{\U_{\theta_{mx}}^-}^{\U^-}\U_{\theta_{mx}}^-)\\=
\Coind_{\U_{\theta_{mx}}^-}^{\U^-}\U_{\theta_{mx}}^-.
\end{gather*}
Now note that the right $\U^-$-module structure on
$\hom_{\U_k}^\c(M(\lambda),
S_{\U_{\theta_{mx}}^-}^{\U_k})$ coincides with the one obtained from
the $\U^-$-bimodule
$\Coind_{\U_{\theta_{mx}}^-}^{\U^-}\U_{\theta_{mx}}^-$ (see Lemma~\ref{a})
with the help of the functor $\Coind_{\U^-}^{\U_k}$. This means that the equality
(*) provides the required isomorphism of the right $\U^-$-modules.

Applying the functor $T_{\theta_{mx}}\circ{}^*$
to both sides of the equality (*)
we obtain the statement of the Lemma.
\qed

Next we apply the functor $\Phi_{mx}$ to the complex $B\bul(\lambda)$.
By definition the obtained complex $\Phi_{mx}(B\bul(\lambda))$ consists of direct sums
of quantum twisted Verma modules. On the other hand by Lemma~\ref{ddd} quantum Verma modules
are $\U_{\theta_{mx}}^-$-free, thus up to a $X\times\Z$-grading shift cohomology
spaces of $\Phi_{mx}(B\bul(\lambda))$ coincide with
$
\left(\Ext_{\U_{\theta_{mx}}^-}(L(\lambda),\U_{\theta_{mx}}^-)\right)^*.
$

\sssn
\Prop        \label{shift}
$\Ext_{\U_{\theta_{mx}}^-}(L(\lambda),\U_{\theta_{mx}}^-)^*
=\delta_{\bullet,\lth(\theta_{mx})}L(\lambda)$.

\subsection{Proof of Proposition~\ref{shift}.}
Note first that $T_{\theta_{mx}}(L(\lambda))=L(\lambda)$ since the character of
$L(\lambda)$ is $W$-invariant and an integrable simple module over $\U$ is completely determined
by its character. Thus up to a $X\times\Z$-grading shift we have
$$
\Ext_{\U_{\theta_{mx}}^-}(L(\lambda),\U_{\theta_{mx}}^-)^*=
\Ext_{T_{\theta_{mx}}(\U_{\theta_{mx}}^-)}(L(\lambda),T_{\theta_{mx}}(\U_{\theta_{mx}}^-))^*.
$$
\sssn
\Lemma
We have  a canonical isomorphism of functors in $\U_{\theta_{mx}}^+\mod$:
$$
\Ext_{\U_{\theta_{mx}}^-}(*,\U_{\theta_{mx}}^-)^*
\til{\map}
\Tor^{\U_{\theta_{mx}}^-}(*,\U_{\theta_{mx}}^{-*}).
$$
\dok
It is sufficient to construct the isomorphism on the level of $\hom$ and $\ten$.
Let $M\in\U_{\theta_{mx}}^-\mod$ and $f\in\hom_{\U_{\theta_{mx}}^-}(M,\U_{\theta_{mx}}^-)$.
Then the morphism
$$
M\ten_ {\U_{\theta_{mx}}^-}
\U_{\theta_{mx}}^{-*}\overset{f\ten\operatorname{Id}}{\map}
\U_{\theta_{mx}}^- \ten_
{\U_{\theta_{mx}}^-} \U_{\theta_{mx}}^{-*} = \U_{\theta_{mx}}^{-*}
\overset{1}{\map}\k
$$
provides a pairing $(f, m\ten\phi),\ m\in M,\phi\in\U_{\theta_{mx}}^{-*}$. One
checks directly that the pairing is perfect for $M=\U_{\theta_{mx}}^-$. Taking
a free resolution of $M\in\U_{\theta_{mx}}^-\mod$ we obtain the statement for an
arbitrary $M$.\qed

By Lemma~\ref{c}(ii) for a graded $T_{\theta_{mx}}(\U_{\theta_{mx}}^-)$-module $M$ with a positive filtration by
finite dimensional $T_{\theta_{mx}}(\U_{\theta_{mx}}^-)$-submodules we have
$\Tor^{T_{\theta_{mx}}(\U_{\theta_{mx}}^-)}(M,T_{\theta_{mx}}(\U_{\theta_{mx}}^{-*})
=\delta_{\bullet,\lth(\theta_{mx})}M$.

Note that the module $L(\lambda)\in\O$, thus in particular, like the proof or
Lemma~\ref{c}(iii), the filtration obtained from the grading satisfies the condition
in~\ref{c}(ii) and by the previous Lemma we have
$$
\Ext_{\U_{\theta_{mx}}^-}(L(\lambda),\U_{\theta_{mx}}^-)^*
=\delta_{\bullet,\lth(\theta_{mx})}L(\lambda).
$$
Now make a $X\times\Z$-grading shift. Proposition~\ref{shift} is proved.\qed

\ssn
\Cor
There exists a complex $\til{B}_{\theta_{mx}}\bul(\lambda)$ in $\M_k$ of the form
\begin{multline*}
\ldots
\map
\underset{{v\in W,}\atop{\lth(v)=m}}{\bigoplus}
M_{\theta_{mx}}(\theta_{mx}v\cdot\lambda)\langle-\hgt(\lambda-\theta_{mx}v\cdot\lambda)\rangle
\map\ldots\\ \map
\underset{{v\in W,}\atop{\lth(v)=1}}{\bigoplus}
M_{\theta_{mx}}(\theta_{mx}v\cdot\lambda)\langle-\hgt(\lambda-\theta_{mx}v\cdot\lambda)\rangle\\
\map
M_{\theta_{mx}}(\theta_{mx}\cdot\lambda)\langle
-\hgt(\lambda-\theta_{mx}\cdot\lambda)\rangle\map 0
\end{multline*}
such that $H^{\ne-\lth(\theta_{mx})}(\til{B}_{\theta_{mx}}\bul(\lambda))=0,\
H^{-\lth(\theta_{mx})}(\til{B}_{\theta_{mx}}\bul(\lambda))=L(\lambda).
$
Here as usual $\langle\cdot\rangle$ denotes the shift of the $\Z$-grading.
\qed

Recall that in~\ref{tlength} we have defined the twisted length function on
the affine Weyl group with the twist $w\in W$ by
$
\lth^w(u)=\lth(w^{-1}u)-\lth(w^{-1}).
$

\sssn
\Cor  \label{twistres}
The complex $\til{B}_{\theta_{mx}}(\lambda)[-\lth(\theta_{mx})]$ can be rewritten as follows:
\begin{multline*}
\ldots
\map
\underset{{v\in W,}\atop{\lth^{\theta_{mx}^{-1}}(v)=n}}{\bigoplus}
M_{\theta_{mx}}(v\cdot\lambda)\langle-\hgt(\lambda-v\cdot\lambda)\rangle
\map\ldots\\ \map
\underset{{v\in W,}\atop{\lth^{\theta_{mx}^{-1}}(v)=0}}{\bigoplus}
M_{\theta_{mx}}(v\cdot\lambda)\langle-\hgt(\lambda-v\cdot\lambda)\rangle
\map\ldots\\
\map
\underset{{v\in W,}\atop{\lth^{\theta_{mx}^{-1}}(v)=-\lth(\theta_{mx})+1}}{\bigoplus}
M_{\theta_{mx}}(v\cdot\lambda)\langle-\hgt(\lambda-v\cdot\lambda)\rangle\\
\map
M_{\theta_{mx}}(\theta_{mx}\cdot\lambda)\langle -\hgt(\lambda-\theta_{mx}\cdot\lambda)\rangle\map 0.\qed
\end{multline*}
\sssn
\Def
We call  the complex $\til{B}_{\theta_{mx}}\bul(\lambda)[-\lth(\theta_{mx})]=:
B_{\theta_{mx}}\bul(\lambda)$ the {\em quantum twisted BGG
resolution} of the module $L(\lambda)$ with the twist $mx$.

We conclude this section with the answer for the semiinfinite homolgy of
$L(\lambda)$ over the algebra $T_{\theta_{mx}}(\U^-)$. Note that
$\U^0$ acts naturally on the semiinfinite $\tor$ spaces  over $T_{\theta_{mx}}(\U^-)$
because $T_{\theta_{mx}}(\U^-)$ is  a normal subalgebra in $T_{\theta_{mx}}(\U^{\le0})$
with the quotient algebra equal to $\U^0$.

\sssn       \label{answer}
\Lemma
We have an equality of $\U^0$-modules
$$
\tor_{\si+n}^{T_{\theta_{mx}}(\U^-)}(\underk,L(\lambda))=
\underset{{v\in W,}\atop{\lth^{\theta_{mx}^{-1}}(v)=n}}{\bigoplus}
\k(v\cdot\lambda).
$$
\dok
The statement follows immediately from Corollary~\ref{twistres}, Lemma~\ref{twistedsem}
and Corollary~\ref{torres}.
\qed

\section{Limit procedure and semiinfinite homology of the infinitely
twisted nilpotent affine quantum group}
Let $\overline{\g}^-\oplus\overline{\h}\oplus\overline{\g}^+$ be the
simple Lie algebra corresponding to the root datum
$(\overline{X},\overline{Y},\ldots)$ of the finite type
$\overline{I}$, let $\hat{\g}:=\g\ten\CC[t,t^{-1}]\oplus\CC K$ be the
affine Kac-Moody algebra corresponding to the root datum
$(X,Y,\ldots)$ of the type $I$. Consider the infinitely twisted
nilpotent subalgebra in $\hat{\g}$
$$
\n^{\si}:=\overline{\g}^-\ten\CC[t,t^{-1}]\oplus\overline{\h}\ten
t^{-1}\CC[t^{-1}] =\underset{\alpha\in
R^{\si-}}{\bigoplus}\hat{\g}_{\alpha}.
$$
In [FF] Feigin and Frenkel constructed semi-infinite BGG resoluitons
of integrable simple modules $L(\lambda)$ over $\hat{\g}$. These
complexes consist  of direct sums of so called Wakimoto modules and
provide semijective resolutions of $L(\lambda)$ over $\n^{\si}$. The
semiinfinite homology spaces of $\n^{\si}$ with coefficients in
$L(\lambda)$ were found this way.

In this section we consider the quantum analogue of the infinitely twisted
nilpotent algebra $\n^{\si}$ and calculate its semiinfiinte homology
with coefficients in integrable simple modules $L(\lambda)$ over $\U$.
Unfortunately we lack the construction of morphisms between quantum
twisted BGG resolutions that would form a porjective system of
complexes like in the affine Lie algebra case (see [Ar3], section 6).
Instead of it using the results of the fifth section we prove that the
semiinfinite $\tor$ spaces themselves form a projective system with
the limit equal to the semiinfinite homology spaces of the infinitely
twisted nilpotent affine quantum group.

\subsection{Infinitely twisted nilpotent affine quantum group.}
Consider a $\k$-vector subspace $\U^{\si-}$ in $\U$ spanned by the
set of PBW monomials as follows:
$$
\left\{\underset{\beta\in\til{R}^{\si-}\cap\til{R}^+}
{\prod}\dot{E}_{\beta}^{a_{\beta}}
\underset{\beta\in\til{R}^{\si-}\cap\til{R}^-}
{\prod}\dot{F}_{\beta}^{b_{\beta}}|
a_{\beta},b_{\beta}\ge0\right\}.
$$
\sssn
\Lemma
$\U^{\si-}$ is a subalgebra in $\U$.

\dok
It is sufficient to check that the product
$\dot{F}_{\beta}\dot{E}_{\alpha}$ of any generators
$\dot{E}_{\alpha},\ \alpha\in\til{R}^{\si-}\cap\til{R}^+$, and
$\dot{F}_{\beta},\ \beta\in\til{R}^{\si-}\cap\til{R}^-$, belongs to
$\U^{\si-}$.  Note that
$$
\til{R}^{\si-}\cap\til{R}^+=\underset{m}{\bigcup}\left((\theta_{mx}(\til{R}^-)\cap\til{R}^+\right)
\text{ and }
\til{R}^{\si-}\cap\til{R}^-=\underset{m}{\bigcap}\left((\theta_{mx}(\til{R}^-)\cap\til{R}^-\right).
$$
On the level of vector spaces we have
\begin{gather*}
\U^{\si-}\cap\U^+=\underset{m}{\bigcup}\left(T_{\theta_{mx}}(\U^-)\cap\U^+\right),\
\U^{\si-}\cap\U^-=\underset{m}{\bigcap}\left(T_{\theta_{mx}}(\U^-)\cap\U^-\right),\\
\text{ and }
\U^{\si-}=\U^{\si-}\cap\U^+\ten\U^{\si-}\cap\U^-.
\end{gather*}
There exists $m_0\ge0$ such that for any integer $m>m_0$ both
$\dot{E}_{\alpha}$ and $\dot{F}_{\beta}$ belong to
$T_{\theta_{mx}}(\U^-)$. Decomposing
$\dot{F}_{\beta}\dot{E}_{\alpha}$ into a sum of tensor product
monomials in $T_{\theta_{mx}}(\U^-)=
T_{\theta_{mx}}(\U^-)\cap\U^+\ten
T_{\theta_{mx}}(\U^-)\cap\U^- $ we obtain the statement of the
Lemma.
\qed

\sssn
\Lemma    \label{dd}
The algebra $\U^{\si-\sharp}$ is isomorphic  to $\U^{\si-}$.

\dok
By~\ref{sind} we have the inclusion of algebras
$\U^{\si-\sharp}\subset\U_k^{\sharp}$ and the image of
$\U^{\si-\sharp}$ can be described as the subalgebra in
$\U_k^{\sharp}$ generated by $\U^{\si-}\cap\U^-$ and
$\U^{\si-}\cap\U^+$. By Theorem~\ref{bimodule} the algebra
$\U_k^{\sharp}$ is isomorphic to $\U_{-2h^{\vee}-k}$. But the
subalgebra generated by $\U^{\si-}\cap\U^-$ and $\U^{\si-}\cap\U^+$
does not depend on the central character and coincides with
$\U^{\si-}$.
\qed

\subsection{Limit procedure.}
Fix integers $m_2>m_1\ge0$. Consider the algebra
$T_{\theta_{m_1x}}(\U^-)\cap
T_{\theta_{m_2x}}(\U^-)$
and two inclusions
\begin{multline*}
i_{m_1}:\ T_{\theta_{m_1x}}(\U^-)\cap
T_{\theta_{m_2x}}(\U^-)
\hookrightarrow
T_{\theta_{m_1x}}(\U^-)\\
\text{ and }
i_{m_2}:\ T_{\theta_{m_1x}}(\U^-)\cap
T_{\theta_{m_2x}}(\U^-)
\hookrightarrow
T_{\theta_{m_2x}}(\U^-).
\end{multline*}
Then
$T_{\theta_{m_1x}}(\U^-)\cap
T_{\theta_{m_2x}}(\U^-)$
has a natural triangular decomposition
$$
T_{\theta_{m_1x}}(\U^-)\cap
T_{\theta_{m_2x}}(\U^-)
=
\U^-\cap
T_{\theta_{m_2x}}(\U^-)
\ten
T_{\theta_{m_1x}}(\U^-)\cap
\U^+.
$$
Note that by the same arguments as in Lemma~\ref{dd} the algebra
$\left(T_{\theta_{m_1x}}(\U^-)\cap
T_{\theta_{m_2x}}(\U^-)\right)^{\sharp}$
coincides with
$T_{\theta_{m_1x}}(\U^-)\cap
T_{\theta_{m_2x}}(\U^-)$.
Next, $i_{m_1}$ is a negative extension of algebras and $i_{m_2}$ is a positive one.
By Propositions~\ref{pos} and~\ref{neg} we obtain canonical morphisms of semiinfinite
$\tor$ spaces:
$$
\Tors^
{T_{\theta_{m_2x}}(\U^-)}
(\underk,L(\lambda))
\map
\Tors^
{T_{\theta_{m_1x}}(\U^-)\cap
T_{\theta_{m_2x}}(\U^-)}
(\underk,L(\lambda))
\map
\Tors^
{T_{\theta_{m_1x}}(\U^-)}
(\underk,L(\lambda)).
$$
Denote the composition map by $p_{m_2,m_1}$.

\sssn
\Prop       \label{projlim}
$
\Tors^
{\U^{\si-}}
(\underk,L(\lambda))
=\invlim
\Tors^
{T_{\theta_{mx}}(\U^-)}
(\underk,L(\lambda)).
$

\dok
First note that by definition of the semiinfinite $\tor$ functor the statement is
equivalent to the following one:
$$
\Exts_
{\U^{\si-}}
(\underk,L(\lambda)^*)
=\dirlim
\Exts_
{T_{\theta_{mx}}(\U^-)}
(\underk,L(\lambda)^*).
$$
Choose  a semijective resolution $M\bul$ of $L(\lambda)^*$ in  $\upC(\U)$  such that
with the differential forgotten it is equal to a $\U$-module of the form
$(S_{\U^+}^{\U})^*\ten V=\Coind_{\U^+}^{\U}\Ind_{\k}^{\U^+}V$. By Proposition~\ref{derived}
we have
$$
\Exts_
{\U^{\si-}}
(\underk,L(\lambda)^*)
=\Exts_
{\U^{\si-}}
(\underk,M\bul)
\text{ and }
\Exts_
{T_{\theta_{mx}}(\U^-)}
(\underk,L(\lambda)^*)=
\Exts_
{T_{\theta_{mx}}(\U^-)}
(\underk,M\bul).
$$
By Lemma~\ref{extres}(iii)
$\Exts_
{\U^{\si-}}
(\underk,M\bul)$
(resp.
$\Exts_
{T_{\theta_{mx}}(\U^-)}
(\underk,M\bul)$)
are the cohomology spaces of the complex
$$
\Hom_
{\U^{\si-}}
(\underk,M\bul\ten_{\U^{\si-}}
S_{\U^{\si-}\cap\U^+}^{\U^{\si-}})
\text{ (resp. }
\Hom_
{T_{\theta_{mx}}(\U^-)}
(\underk,M\bul\ten_{T_{\theta_{mx}}(\U^-)}
S_{T_{\theta_{mx}}(\U^-)\cap\U^+}^{T_{\theta_{mx}}(\U^-)})).
$$
Next note that
\begin{gather*}
\Hom_
{\U^{\si-}}
(\underk,(S_{\U^+}^{\U})^*\ten_{\U^{\si-}}
S_{\U^{\si-}\cap\U^+}^{\U^{\si-}})
=\U^{\si+}\cap\U^+\ten(\U^{\si+}\cap\U^-)^*
\text{ and }\\
\Hom_
{T_{\theta_{mx}}(\U^-)}
(\underk,(S_{\U^+}^{\U})^*\ten_{T_{\theta_{mx}}(\U^-)}
S_{T_{\theta_{mx}}(\U^-)\cap\U^+}^{T_{\theta_{mx}}(\U^-)})=
T_{\theta_{mx}}(\U^+)\cap\U^+
\ten
(T_{\theta_{mx}}(\U^+)\cap\U^-)^*.
\end{gather*}
Note also that the canonical semiinfiinte $\ext$ morphisms here coincide with the obvious ones.
In particular we have
$$
\Exts_ {\U^{\si-}}
(\underk,(S_{\U^+}^{\U})^*)= \dirlim
\Exts_ {T_{\theta_{mx}}(\U^-)}
(\underk,(S_{\U^+}^{\U})^*).
$$
Thus the Proposition is proved since the $\dirlim$ functor is exact.
\qed

From Lemma~\ref{answer} we know that the space $\tor_{\si+n}^{T_{\theta_{mx}}(\U^-)}
(\underk,L(\lambda))$ has  a base enumerated by the elements of the affine Weyl group
$w\in W$ such that $\lth^{\theta_{mx}^{-1}}(w)=n$. Next recall that by Corollary~\ref{qq}
for every $w\in W$ there exists $m_0\in\N$ such that for every $m>m_0$ we have
$\lth^{\theta_{mx}^{-1}}(w)=\lth^{\si}(w)$. Fix $w\in W$ and choose $m>m_0=m_0(w)$.
Let $n=\lth^{\si}(w)$. Consider the map
\begin{multline*}
p_{m+1,m}:\
\tor^
{T_{\theta_{(m+1)x}}(\U^-)}_{\si+n}
(\underk,L(\lambda))
\map
\tor^
{T_{\theta_{mx}}(\U^-)\cap
T_{\theta_{(m+1)x}}(\U^-)}_{\si+n}
(\underk,L(\lambda))\\
\map
\tor^
{T_{\theta_{mx}}(\U^-)}_{\si+n}
(\underk,L(\lambda)).
\end{multline*}
Denote the first (resp. the second) map in the composition
by $p_{m+1,m}^-$ (resp. by $p_{m+1,m}^+$).
Denote the base vector corresponding to $w$ in
$\tor^
{T_{\theta_{(m+1)x}}(\U^-)}_{\si+n}
(\underk,L(\lambda))$
(resp. in
$\tor^
{T_{\theta_{mx}}(\U^-)}_{\si+n}
(\underk,L(\lambda))$)
by $a_{m+1}^w$ (resp. by $a_m^w$).

\sssn
\Theorem \label{stabil}
$p_{m+1,m}(a_{m+1}^w)=c a_m^w$ and $c\ne0$ for $m>>0$.

\dok
Since the map $p_{m+1,m}$ preserves the natural $X$-gradings on the semiinfinite $\tor$ spaces,
we have to prove only that $p_{m+1,m}(a_{m+1}^w)\ne0$ for $m>>0$.
Note that $B_{\theta_{mx}}\bul(\lambda)$ (resp. $B_{\theta_{(m+1)x}}\bul(\lambda)$)
is a co-semijective resolution of $L(\lambda)$ not only over
$T_{\theta_{mx}}(\U^-)$
(resp. over
$T_{\theta_{(m+1)x}}(\U^-)$)
but also over
$T_{\theta_{mx}}(\U^-)\cap
T_{\theta_{(m+1)x}}(\U^-)$.
Thus by Proposition~\ref{derived} we have
\begin{multline*}
\tor^
{T_{\theta_{mx}}(\U^-)\cap
T_{\theta_{(m+1)x}}(\U^-)}_{\si+n}
(\underk,L(\lambda))=
\tor^
{T_{\theta_{mx}}(\U^-)\cap
T_{\theta_{(m+1)x}}(\U^-)}_{\si+n}
(\underk,B_{\theta_{mx}}\bul(\lambda))\\ =
\tor^
{T_{\theta_{mx}}(\U^-)\cap
T_{\theta_{(m+1)x}}(\U^-)}_{\si+n}
(\underk,B_{\theta_{(m+1)x}}\bul(\lambda)).
\end{multline*}
Consider the morphisms of the complexes
\begin{gather*}
p_{m+1,m}^-:\
\Tors^
{T_{\theta_{(m+1)x}}(\U^-)
}
(\underk,B_{\theta_{(m+1)x}}(\lambda))
\map
\Tors^
{T_{\theta_{mx}}(\U^-)\cap
T_{\theta_{(m+1)x}}(\U^-)}
(\underk,B_{\theta_{(m+1)x}}(\lambda))\\
\text{ and }
\Tors^
{T_{\theta_{mx}}(\U^-)\cap
T_{\theta_{(m+1)x}}(\U^-)}
(\underk,B_{\theta_{mx}}(\lambda))
\map
\Tors^
{T_{\theta_{mx}}(\U^-)
}
(\underk,B_{\theta_{mx}}(\lambda)),
\end{gather*}
the differentials in the complexes obtained from the ones in
$B_{\theta_{(m+1)x}}\bul(\lambda)$ and $B_{\theta_{mx}}\bul(\lambda)$ respectively.
Note that the differentials in
$$
\Tors^
{T_{\theta_{mx}}(\U^-)}
(\underk,B_{\theta_{mx}}(\lambda))
\text{ and }
\Tors^
{T_{\theta_{(m+1)x}}(\U^-)
}
(\underk,B_{\theta_{(m+1)x}}(\lambda))
$$
vanish.
By Lemma~\ref{canmap2} up to a $X$-grading shifts the map
$$
p_{m+1,m}^+:\
\tor_{\si+0}^
{T_{\theta_{mx}}(\U^-)\cap
T_{\theta_{(m+1)x}}(\U^-)}
(\underk,B_{\theta_{mx}}^n(\lambda))
\map
\tor_{\si+0}^
{T_{\theta_{mx}}(\U^-)
}
(\underk,B_{\theta_{mx}}^n(\lambda))
$$
coincides with the canonical map
$$
\underset{{v\in W,}\atop{\lth^{\theta_{mx}^{-1}}(v)=n}}{\bigoplus}
\Ind_
{T_{\theta_{(m+1)x}}(\U^-)\cap\U^-}^
{T_{\theta_{mx}}(\U^-)\cap\U^-}
\underk
\map
\underset{{v\in W,}\atop{\lth^{\theta_{mx}^{-1}}(v)=n}}{\bigoplus}
\underk.
$$
In particular it is surjective.
By Lemma~\ref{canmap1} up to a $X$-grading shifts the map
$$
p_{m+1,m}^-:\
\tor_{\si+0}^
{T_{\theta_{(m+1)x}}(\U^-)
}
(\underk,B_{\theta_{(m+1)x}}^n(\lambda))
\map
\tor_{\si+0}^
{T_{\theta_{mx}}(\U^-)\cap
T_{\theta_{(m+1)x}}(\U^-)}
(\underk,B_{\theta_{(m+1)x}}^n(\lambda))
$$
coincides with the canonical map
$$
\underset{{v\in W,}\atop{\lth^{\theta_{(m+1)x}^{-1}}(v)=n}}{\bigoplus}\underk
\map
\underset{{v\in W,}\atop{\lth^{\theta_{(m+1)x}^{-1}}(v)=n}}{\bigoplus}
\Coind_{\U_{\theta_{mx}}^+}^{\U_{\theta_{(m+1)x}}^+}\underk.
$$
Our next aim is to prove that $p_{m+1,m}(a_{m+1}^w)$ defines
a nonzero cohomology class.
Denote the component of the differential $M_{\theta_{mx}}(v\cdot\lambda)\map
M_{\theta_{mx}}(w\cdot\lambda)$ in $B_{\theta_{mx}}\bul(\lambda)$ by $d_{v,w}$.

\sssn
\Lemma
For every $w\in W$ there exists $m_1\in \N$ such that for every $m>m_1$ the set
$
Prev(w,\lambda,m):=\{v\in W|\lth^{\theta_{mx}^{-1}}(v)=\lth^{\theta_{mx}^{-1}}(w)-1,
d_{v,w}\ne0\}
$
(resp. the set
$
Next(w,\lambda,m):=\{v\in W|\lth^{\theta_{mx}^{-1}}(v)=\lth^{\theta_{mx}^{-1}}(w)+1,
d_{w,v}\ne0\})
$
coincides with the set
$
\{v\in W|\lth^{\si}(v)=\lth^{\si}(w)-1,
v\le^{\si}w\}
\text{ (resp. }
\{v\in W|\lth^{\si}(v)=\lth^{\si}(w)+1,
w\le^{\si}v\}).
$
\qed

It is known that for any $w\in W$ both sets from the previous Lemma are finite.
We suppose that $m>m_1=m_1(w)$. Let
$$
\operatorname{ht}^-:=\underset{v\in Prev(w,\lambda,m)}{\operatorname{max}}
|\hgt(w\cdot\lambda-v\cdot\lambda)|,\
\operatorname{ht}^+:=\underset{v\in Next(w,\lambda,m)}{\operatorname{max}}
|\hgt(w\cdot\lambda-v\cdot\lambda)|.
$$
Note that the space
$\Coind_{\U_{\theta_{mx}}^+}^{\U_{\theta_{(m+1)x}}^+}\underk$
has a base of dual vectors to PBW monomials with generators $\dot{E}_{\alpha},\
\alpha\in R_{\theta_{(m+1)x}}^+\setminus R_{\theta_{mx}}^+$.
Choosing $m$ big enough we obtain $\hgt(\alpha)>\hgt^-$ for every
$\alpha\in R_{\theta_{(m+1)x}}^+\setminus R_{\theta_{mx}}^+$.
Thus for $m$ big enough the vector
$p_{m+1,m}^-(a_{m+1}^w)$ provides a nontrivial class in the cohomology of the complex
$\tor_{\si+0}^
{T_{\theta_{mx}}(\U^-)\cap
T_{\theta_{(m+1)x}}(\U^-)}
(\underk,B_{\theta_{(m+1)x}}(\lambda))$. It is proved similarly
that $a_m^w$ belongs to
$$
\im(\tor^{T_{\theta_{mx}}(\U^-)\cap
T_{\theta_{(m+1)x}}(\U^-)}_{\si+n}(\underk,L(\lambda))
\map\tor^{T_{\theta_{mx}}(\U^-)}_{\si+n}(\underk,L(\lambda)))
$$
Now note that the grading component in
$\Coind_{\U_{\theta_{mx}}^+}^{\U_{\theta_{(m+1)x}}^+}\k(w\cdot\lambda)$
of the weight $w\cdot\lambda$ is one dimensional. Thus $p_{m+1,m}$ nessesarily
takes $a_{m+1}^w$ to a nonzero vector $ca_m^w$.
\qed

\sssn
\Cor        \label{final}
We have an equality of $\U^0$-modules
$$
\tor_{\si+n}^{\U^{\si-}}(\underk,L(\lambda))=
\underset{{v\in W,}\atop{\lth^{\si}(v)=n}}{\bigoplus}
\k(v\cdot\lambda).
$$
\dok
The statement follows immediately from Lemma~\ref{answer}, Proposition~\ref{projlim}
and Theorem~\ref{stabil}.
\qed

\section*{References}
$\text{[Ar1]}$ S.~Arkhipov. {\it Semiinfinite cohomology of quantum groups.}
Preprint q-alg/9601026 (1996), 1-24.\\
$\text{[Ar2]}$ S.~Arkhipov. {\it
Semiinfinite cohomology of associative algebras and bar duality.} Preprint
q-alg/9602013 (1996), 1-21.\\
$\text{[Ar3]}$ S.~Arkhipov. {\it
A new construction of the semi-infinite BGG resolution.} Preprint
q-alg/9605043 (1996), 1-29.\\
$\text{[Be]}$ J. Beck. {\em Convex bases of PBW type for quantum affine algebras},
Comm. Math. Phys. \hbox{\bf 165} (1994), 193-199.\\
$\text{[BeK]}$ J. Beck, V. Kac. {\em Finite dimensional representations of
quantum affine algebras at roots of unity.} Preprint (1996), 1-32.\\
$\text{[BGG]}$ I.~I.~Bernstein, I.~M.~Gelfand,
S.~I.~Gelfand.  {\it Differential operators on the principal affine space and
investigation of} $\g${\em -modules}, in Proceedings of Summer School on Lie
groups of Bolyai J\=anos Math. Soc., Helstead, New York, 1975.\\
$\text{[DCKP] C.De Concini, V.Kac, C.Procesi.}$ {\it Some
remarkable degenerations of quantum groups.} Comm. Math. Phys.
\hbox{\bf 157}  (1993), p.405-427.\\
$\text{[F]}$ B. Feigin.
{\it Semi-infinite cohomology of Kac-Moody and Virasoro Lie algebras.} Usp.
Mat.  Nauk \hbox{\bf 33}, no.2 (1984),  195-196 (in Russian).\\
$\text{[FF]}$ B. Feigin, E. Frenkel. {\em Affine Kac-Moody algebras and
semi-infinite flag manifolds.} Comm. Math. Phys. \hbox{\bf 128} (1990),
617-639.\\
$\text{[K]}$ V. Kac. {\it Infinite dimensional Lie algebras.}
Birkh\"auser, Boston (1984).\\
$\text{[L1]}$ G.~Lusztig. {\em Introduction to quantum groups.} (Progress in
Mathematics \hbox{\bf 110}), Boston etc. 1993 (Birkh\"auser).\\
$\text{[L2]}$ G. Lusztig. {\it Hecke algebras and Jantzen's
generic decomposition patterns.} Adv. in Math. Vol.\hbox{\bf 37},
No.\hbox{\bf 2} (1980), 121-164.\\
$\text{[M]}$ F. Malikov. {\em Quantum groups: singular vectors and BGG resolution.}
KURIMS preprint (1991), 1-21.\\
$\text{[P]}$ P. Papi. {\em Convex orderings in affine root systems}. To appear in
Journal of Algebra.\\
$\text{[Pr]}$ S. Priddy. {\em Koszul resolutions.} Trans. AMS \hbox{\bf 152}, no.1
(1970), 39-60.\\
A.~Voronov. {\it Semi-infinite
homological algebra.} Invent. Math. {\bf 113}, (1993),  103-146.

\end{document}